\documentclass[11pt]{article}
\usepackage{epic,latexsym,amssymb}
\usepackage{color}
\usepackage{tikz}
\usepackage{amsmath}
\usepackage{comment}
\usepackage{geometry,graphicx}

\newtheorem{theorem}{Theorem}[section]
\newtheorem{lemma}[theorem]{Lemma}
\newtheorem{quest}[theorem]{Question}

\newtheorem{claim}{Claim}
\newtheorem{subclaim}{Claim}[claim]

\newtheorem{corollary}[theorem]{Corollary}

\newtheorem{conj}[theorem]{Conjecture}

\newcommand{\QED}{$\Box$}

\newcommand{\smallqed}{{\tiny ($\Box$)}}

\newcommand{\dom}{{\rm dom}}

\newcommand{\pd}{{\rm pd}}

\newcommand{\pn}{\mbox{pn}}

\newcommand{\epn}{{\rm epn}}

\newcommand{\proof}{\noindent\textbf{Proof. }}

\newcommand{\2}{ \vspace{0.2cm} }
\newcommand{\1}{ \vspace{0.1cm} }

\let\oldenumerate\enumerate
\renewcommand{\enumerate}{
  \oldenumerate
  \setlength{\itemsep}{0pt}
  \setlength{\parskip}{0pt}
  \setlength{\parsep}{0pt}
}

\def\vertex(#1){\put(#1){\circle*{2}}}
\def\vertexo(#1){\put(#1){\circle{2}}}
\def\vert(#1){\put(#1){\circle*{1.5}}}
\def\verto(#1){\put(#1){\circle{1.5}}}
\def\lab(#1)#2{\put(#1){\makebox(0,0)[c]{#2}}}

\begin{document}

\title{Partial domination in supercubic graphs}

\author{$^{1,2}$Csilla Bujt\'{a}s \qquad $^3$Michael A. Henning \qquad
$^{1,2,4}$Sandi Klav\v zar
\\ \\
\small $^1$Faculty of Mathematics and Physics \\
\small University of Ljubljana, Slovenia\\
\small \tt csilla.bujtas@fmf.uni-lj.si\\
\small \tt sandi.klavzar@fmf.uni-lj.si\\
\\
\small $^2$Institute of Mathematics, Physics and Mechanics \\
\small Ljubljana, Slovenia\\
\\
\small $^3$Department of Mathematics and Applied Mathematics \\
\small University of Johannesburg \\
\small Auckland Park, 2006 South Africa\\
\small \tt mahenning@uj.ac.za\\
\\
\small $^4$ Faculty of Natural Sciences and Mathematics\\
\small University of Maribor, Slovenia
}
\date{}
\maketitle

\begin{abstract}
For some $\alpha$ with $0 < \alpha \le 1$, a subset $X$ of vertices in a graph $G$ of order~$n$ is an $\alpha$-partial dominating set of $G$ if the set $X$ dominates at least $\alpha \times n$ vertices in $G$. The $\alpha$-partial domination number ${\rm pd}_{\alpha}(G)$ of $G$ is the minimum cardinality of an $\alpha$-partial dominating set of $G$. In this paper partial domination of graphs with minimum degree at least $3$ is studied. It is proved that if $G$ is a graph of order~$n$ and with $\delta(G)\ge 3$, then ${\rm pd}_{\frac{7}{8}}(G) \le \frac{1}{3}n$. If in addition  $n\ge 60$, then ${\rm pd}_{\frac{9}{10}}(G) \le \frac{1}{3}n$, and if $G$ is a connected cubic graph of order $n\ge 28$, then ${\rm pd}_{\frac{13}{14}}(G) \le \frac{1}{3}n$. Along the way it is shown that there are exactly four connected cubic graphs of order $14$ with domination number $5$.
\end{abstract}


\noindent
\textbf{AMS Subj.\ Class.:} 05C69

\newpage
\section{Introduction}

One of the central themes of the theory of graph domination is establishing upper bounds for graphs with a prescribed minimum degree as a function of graph order. The topic is in depth surveyed in the paper~\cite{Henning-2022} as well as in the 2023 book~\cite{HaHeHe-23}. Special attention has been paid to cubic graphs and graphs of minimum degree at least $3$. For the latter graphs, Reed~\cite{Re-96} in 1996 established a best possible upper bound.

\begin{theorem}{\rm (\cite{Re-96})}
\label{thm:Reed}
If $G$ is a graph of order $n$ with $\delta(G) \ge 3$, then $\gamma(G) \le \frac{3}{8}n$.
\end{theorem}

In 2009, Kostochka and Stocker sharpened Reed's bound for connected cubic graphs as follows.

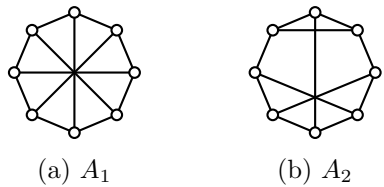
\begin{figure}[ht!]
\begin{center}
\begin{tikzpicture}[scale=.8,style=thick,x=0.8cm,y=0.8cm]
\def\vr{2.5pt} 
\path (0.00,1.25) coordinate (u1);
\path (0.37,0.37) coordinate (u2);
\path (0.37,2.13) coordinate (u8);
\path (1.25,0.00) coordinate (u3);
\path (1.25,2.50) coordinate (u7);
\path (2.13,0.37) coordinate (u4);
\path (2.13,2.13) coordinate (u6);
\path (2.50,1.25) coordinate (u5);
%
\draw (u1)--(u2)--(u3)--(u4)--(u5)--(u6)--(u7)--(u8)--(u1);
\draw (u1)--(u5);
\draw (u2)--(u6);
\draw (u3)--(u7);
\draw (u4)--(u8);
\draw (u1) [fill=white] circle (\vr);
\draw (u2) [fill=white] circle (\vr);
\draw (u3) [fill=white] circle (\vr);
\draw (u4) [fill=white] circle (\vr);
\draw (u5) [fill=white] circle (\vr);
\draw (u6) [fill=white] circle (\vr);
\draw (u7) [fill=white] circle (\vr);
\draw (u8) [fill=white] circle (\vr);
\draw (1.25,-0.8) node {{\small (a) $A_1$}};
\path (5.00,1.25) coordinate (v1);
\path (5.37,0.37) coordinate (v2);
\path (5.37,2.13) coordinate (v8);
\path (6.25,0.00) coordinate (v3);
\path (6.25,2.50) coordinate (v7);
\path (7.13,0.37) coordinate (v4);
\path (7.13,2.13) coordinate (v6);
\path (7.50,1.25) coordinate (v5);
%
\draw (v1)--(v2)--(v3)--(v4)--(v5)--(v6)--(v7)--(v8)--(v1);
\draw (v6)--(v8);
\draw (v3)--(v7);
\draw (v1)--(v4);
\draw (v2)--(v5);
\draw (v1) [fill=white] circle (\vr);
\draw (v2) [fill=white] circle (\vr);
\draw (v3) [fill=white] circle (\vr);
\draw (v4) [fill=white] circle (\vr);
\draw (v5) [fill=white] circle (\vr);
\draw (v6) [fill=white] circle (\vr);
\draw (v7) [fill=white] circle (\vr);
\draw (v8) [fill=white] circle (\vr);
\draw (6.25,-0.8) node {{\small (b) $A_2$}};
\end{tikzpicture}
\caption{The (non-planar) cubic graphs $A_1$ and $A_2$ of order $n = 8$ with $\gamma(A_1) = \gamma(A_2) = 3$}
\label{f:cubic8}
\end{center}
\end{figure}

\begin{theorem}{\rm (\cite{KoStocker-09})}
\label{thm:KoSt}
If $G$ is a connected cubic graph of order~$n$, then $\gamma(G) \le \frac{5}{14}n$, unless $G$ is one of the two graphs $A_1$ and $A_2$ shown in Figure~\ref{f:cubic8}.
\end{theorem}

Kostochka and Stocker further proved that the graphs $A_1$ and $A_2$ are the only connected, cubic graphs that achieve the $\frac{3}{8}$-bound of Theorem~\ref{thm:Reed}. On the other hand, Reed~\cite{Re-96} conjectured that $\gamma(G) \le \lceil \frac{1}{3}n \rceil$ whenever $G$ is a connected cubic graph of order~$n$. Kostochka and Stodolsky~\cite{KoSt-05} disproved this conjecture by constructing an infinite sequence $\{G_k\}^\infty_{k=1}$ of connected, cubic graphs with
\[
\lim_{k\to\infty} \frac{\gamma(G_k)}{|V(G_k)|} \ge \frac{1}{3} + \frac{1}{69}.
\]
Subsequently, Kelmans~\cite{Ke-06} constructed an infinite series of $2$-connected, cubic graphs $H_k$ with
\[
\lim_{k\to\infty} \frac{\gamma(H_k)}{|V(H_k)|} \ge \frac{1}{3} + \frac{1}{60}.
\]
Thus, there exist connected cubic graphs $G$ of arbitrarily large order~$n$ satisfying
\[
\gamma(G) \ge \left( \frac{1}{3} + \frac{1}{60} \right)n.
\]

So, $\gamma(G) \le \lceil \frac{1}{3}n \rceil$ does not hold for all connected cubic graphs. On the other hand, in 2010 Verstra\"{e}te conjectured that if $G$ is a cubic graph of order $n$ and girth at least $6$, then $\gamma(G) \le \frac{1}{3}n$, see~\cite[Conjecture~10.23]{HaHeHe-23}. In~\cite{Lowenstein-2008} the conjecture has been verified for cubic graphs with girth at least $83$. Further upper bounds on the domination number of a cubic graph in terms of its order and girth were proved in~\cite{KoSt-09, Kral-2012, Rautenbach-2008}.

The following concepts were independently introduced in~\cite{Case-2017, Das-2019}. Let $G = (V(G), E(G))$ be a graph of order~$n$. For some $\alpha$ with $0 < \alpha \le 1$, a set $S\subseteq V(G)$ is an $\alpha$-\emph{partial dominating set} of $G$ if
\begin{equation*}
\label{Eq1}
|N_G[S]| \ge \alpha \times n, \1
\end{equation*}
that is, the set $S$ dominates at least $\alpha n$ vertices in $G$.
The \emph{$\alpha$-partial domination number} of $G$, denoted by $\pd_{\alpha}(G)$ (also by $\gamma_{\alpha}(G)$ in the literature), is the minimum cardinality of an $\alpha$-partial dominating set of $G$. Investigations on the concept of partial domination in graphs can be found  in~\cite{CaIs-21,Case-2017,Das-2019,DaDe-2018,MaIs-10,PhVa-21}. At this point, it should be pointed out that the term ``partial domination" is also used to refer to a concept that is different from ours~\cite{Borg-2020}. We also remark that the concept of an $\alpha$-dominating set~\cite{DaLa-2018,Dunbar-2000,RadVolk-2016} is different from our concept of an $\alpha$-partial dominating set.

In light of the above, this paper addresses the following natural question: What is the largest possible value on $\alpha$ such that the $\alpha$-partial domination number of a connected cubic graph is at most one-third the order of the graph? We further consider the same question in the more general setting of graphs with minimum degree at least~$3$.

We proceed as follows. In Section~\ref{S:prelim}, we present the graph theory terminology we adopt in this paper, and state preliminary results. In Section~\ref{sec:order4}, we prove that the $\frac{7}{8}$-partial domination number of a connected cubic graph $G$ of order~$14$ is at most~$4$. Thereafter in Section~\ref{sec:delta-ge-3}, we prove that the $\frac{7}{8}$-partial domination number of a graph with minimum degree at least~$3$ is at most one-third its order, and prove a stronger statement if the order of the graph is large enough. In Section~\ref{S:closing}  we show that there are exactly four connected cubic graphs of order~$14$ with domination number~$5$, and conjecture that these are the only graphs achieving equality in the upper bound $\gamma(G) \le \frac{5}{14}n$ given by Kostochka and Stocker in Theorem~\ref{thm:KoSt}.

\section{Preliminaries}
\label{S:prelim}

In this section, we call up the definitions, concepts and known results that we need for what follows. Let $G = (V(G), E(G))$ be a graph. The \emph{open neighborhood} $N_G(v)$ of a vertex $v$ in $G$ is the set of vertices adjacent to $v$, while the \emph{closed neighborhood} of $v$ is the set $N_G[v] = \{v\} \cup N_G(v)$. For a set $D \subseteq V(G)$, its \emph{open neighborhood} is the set $N_G(D) = \cup_{v \in D} N_G(v)$, and its \emph{closed neighborhood} is the set $N_G[D] = N_G(D) \cup D$. The minimum and maximum degrees in $G$ are denoted by $\delta(G)$ and $\Delta(G)$, respectively. The graph $G$ is $r$-\emph{regular} if every vertex in $G$ has degree~$r$. A $3$-regular graph is called a \emph{cubic graph}, and a graph $G$ with $\Delta(G)\le 3$ a \emph{subcubic graph}. To these established terms we add the term {\em supercubic graph} which refers to graphs $G$ with $\delta(G) \ge 3$.

A \emph{dominating set} of a graph $G$ is a set $S$ of vertices of $G$ such that every vertex not in $S$ has a neighbor in $S$. The \emph{domination number} of $G$, denoted by $\gamma(G)$, is the minimum cardinality of a dominating set. A $\gamma$-set of $G$ is a dominating set of $G$ of minimum cardinality~$\gamma(G)$. Let $X$ and $Y$ be subsets of vertices in $G$. The set $X$ \emph{dominates} the set $Y$ if every vertex in $Y$ is in the set $X$ or has a neighbor in the set $X$, that is, if $Y \subseteq N_G[X]$. If $X$ is a set of vertices in a graph $G$, then we denote by $\dom_G(X)$ the number of vertices dominated by the set $X$, and so $\dom_G(X) = |N_G[X]|$. A thorough treatise on domination in graphs can be found in~\cite{HaHeHe-20,HaHeHe-21,HaHeHe-23}.

For a set of vertices $S$ in a graph $G$ and a vertex $v \in S$, the \emph{$S$-private neighborhood} of $v$ is defined by $\pn[v,S] = \{w \in V(G) \colon N_G[w] \cap S = \{v\}\}$. The \emph{$S$-external private neighborhood} of $v$ is the set $\epn[v,S] = \pn[v,S] \setminus S$. (The set $\epn[v,S]$ is also denoted $\epn(v,S)$ in the literature.) An \emph{$S$-external private neighbor} of $v$ is a vertex in $\epn[v,S]$. In 1979, Bollob\'{a}s and Cockayne~\cite{BoCo-79} established the following property of minimum dominating sets in graphs to be used later on.

\begin{lemma}{\rm (\cite{BoCo-79})}
\label{lem:BoCo}
Every isolate-free graph $G$ contains a $\gamma$-set $D$ such that $\epn[v,D] \ne \emptyset$ for every vertex $v \in D$.
\end{lemma}

A set $S$ of vertices in $G$ is a \emph{packing} in $G$ if the closed neighborhoods of vertices in $S$ are pairwise disjoint. Equivalently, $S$ is a packing in $G$ if the vertices in $S$ are pairwise at distance at least~$3$. A packing is sometimes called a $2$-packing in the literature. The \emph{packing number} of $G$, denoted by $\rho(G)$, is the maximum cardinality of a packing in $G$. In 1996, Favaron~\cite{Fa-96} proved the following result on the packing number of a cubic graph.

\begin{theorem}{\rm (Favaron~\cite{Fa-96})}
\label{t:pack1}
If $G$ is a connected cubic graph of order $n$ different from the Petersen graph, then $\rho(G) \ge \frac{n}{8}$.
\end{theorem}

For a set of vertices in $G$, the subgraph of $G$ induced by $S$ is denoted by $G[S]$. Finally, the \emph{boundary} of a set $S$ of vertices in $G$, denoted by $\partial(S)$, is the set of vertices not in $S$ that have a neighbor in $S$, that is, $\partial(S) = N_G[S] \setminus S$.

\section{(Partial) domination in cubic graphs of order $14$}
\label{sec:order4}

In this section, we present a preliminary result that the $\frac{7}{8}$-partial domination number of a connected cubic graph $G$ of order~$14$ is at most~$4$. We will need this result when proving our main theorem in Section~\ref{sec:delta-ge-3}.

\begin{theorem}
\label{t:cubicn14}
If $G$ is a connected cubic graph of order $n = 14$, then
$$\pd_{\frac{7}{8}}(G) \le 4 < \frac{1}{3}n\,.$$
\end{theorem}
\proof Let $\alpha = \frac{7}{8}$ and let $G$ be a connected cubic graph of order~$n = 14$.  Let $\gamma = \gamma(G)$. By Theorem~\ref{thm:KoSt}, $\gamma \le \lfloor \frac{5}{14}n \rfloor = 5$. If $\gamma \le \lfloor \frac{1}{3}n \rfloor = 4$, then every $\gamma$-set of $G$ is certainly an $\alpha$-partial dominating set of $G$. Thus in this case, $\pd_{\alpha}(G) \le \gamma \le 4$, as desired. Hence we may assume in what follows that $\gamma = 5$.

By Theorem~\ref{t:pack1}, the graph $G$ has packing number $\rho(G) \ge \lceil \frac{n}{8} \rceil = 2$. Let $P$ be a maximum packing in $G$, and so $|P| = \rho(G) \ge 2$. Suppose that $\rho(G) > 2$, implying that $\rho(G) = 3$. In this case, $\dom_G(P) = |N_G[P]| = 12$. Thus if $v$ is any one of the two vertices in $V(G) \setminus N_G[P]$ and $S = P \cup \{v\}$, then the set $S$ satisfies $\dom_G(S) \ge 13 > \frac{7}{8}n$, and so $\pd_{\alpha}(G) \le |S| = 4$. Hence we may assume that $|P| = \rho(G) = 2$, for otherwise the desired result follows.

Let $X = V(G) \setminus N_G[P]$, and so $|X| = 6$. If a vertex in $X$ has all three of its neighbors in the set $X$, then we can add such a vertex to the set $P$ to produce a packing of cardinality~$3$, contradicting our assumption that $\rho(G) = 2$. Hence, every vertex in $X$ has at most two neighbors that belong to $X$.

Suppose that a vertex $x \in X$ has two neighbors in the set $X$. In this case, we let $P_x = P \cup \{x\}$. The resulting set $P_x$ satisfies $|P_x| = 3$ and $\dom_G(P_x) = 4 + 4 + 3 = 11$. Let $Z = V(G) \setminus N_G[P_x]$, and so $|Z| = 3$. If there is a vertex $z \in Z$ with at least one neighbor in $Z$, then the set $S = P_x \cup \{z\}$ satisfies $\dom_G(S) \ge 13$ and $|S| = 4$, and so as before $\pd_{\alpha}(G) \le |S| = 4$. Hence, we may assume that $Z$ is an independent set in $G$. Thus, every vertex in $Z$ has all three of its neighbors contained in the boundary $\partial(P_x)$ of the set $P_x$. Denoting by $\ell_1$ the number of edges between the sets $Z$ and $\partial(P_x)$, we obtain $\ell_1 = 3|Z| = 9$. Since $|\partial(P_x)| = \dom_G(P_x) - |P_x| = 11 - 3 = 8$, by the Pigeonhole Principle at least one vertex $v$ in the boundary $\partial(P_x)$ of $P_x$ has at least two neighbors in $Z$. Thus the set $S = P_x \cup \{v\}$ satisfies $\dom_G(S) \ge 13$ and $|S| = 4$, and so as before $\pd_{\alpha}(G) \le |S| = 4$.

Hence, we may assume that every vertex in $X$ has at most one neighbor that belongs to $X$, and therefore at least two neighbors that belong to the boundary $\partial(P)$ of $P$. Denoting by $\ell_2$ the number of edges between the sets $X$ and $\partial(P)$, we obtain $\ell_2 \ge 2|X| = 2 \times 6 = 12$. However every vertex in $\partial(P)$ has one neighbor in $P$ and therefore at most two neighbors in $X$, and so $\ell_2 \le 2|\partial(P)| = 2 \times 6 = 12$. Consequently, $\ell_2 = 12$, implying that $\partial(P)$ is an independent set and each vertex in $\partial(P)$ has exactly two neighbors in $X$. Furthermore, each vertex in $X$ has exactly two neighbors in $\partial(P)$ and one neighbor in $X$. In particular, the subgraph induced by the set $X$ consists of three disjoint copies of $P_2$, that is, $G[X] = 3P_2$.

Let $Y = \partial(P)$, and let $H$ be the graph with vertex set $X \cup Y$ and with edge set consisting of all edges in $G$ between $X$ and $Y$. By our earlier observations, $|X| = |Y| = 6$. The resulting bipartite graph $H$ has partite sets $X$ and $Y$ and is a $2$-regular graph of order~$12$. Thus, either $H = 2C_6$, or $H = 3C_4$, or $H = C_4 \cup C_8$, or $H = C_{12}$. Let $P = \{v_1,v_2\}$. Let $X = \{x_1,x_2,\ldots,x_6\}$ and $Y = \{y_1,y_2,\ldots,y_6\}$.

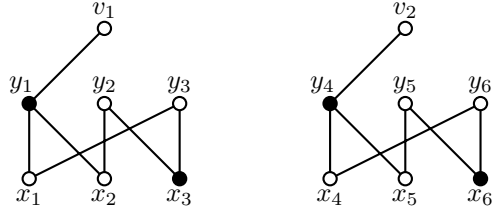
\begin{figure}[ht!]
\begin{center}
\begin{tikzpicture}[scale=1,style=thick,x=1cm,y=1cm]
\def\vr{2.5pt} 
\path (0,0) coordinate (x1);
\path (0,1) coordinate (y1);
\path (-0.1,1) coordinate (y1p);
\path (1,0) coordinate (x2);
\path (1,1) coordinate (y2);
\path (2,0) coordinate (x3);
\path (2,1) coordinate (y3);
\path (1,2) coordinate (v1);
\path (4,0) coordinate (x4);
\path (4,1) coordinate (y4);
\path (3.9,1) coordinate (y4p);
\path (5,0) coordinate (x5);
\path (5,1) coordinate (y5);
\path (6,0) coordinate (x6);
\path (6,1) coordinate (y6);
\path (5,2) coordinate (v2);
\draw (x1)--(y1)--(x2)--(y2)--(x3)--(y3)--(x1);
\draw (v1)--(y1);
\draw (x4)--(y4)--(x5)--(y5)--(x6)--(y6)--(x4);
\draw (v2)--(y4);
%
\draw (v1) [fill=white] circle (\vr);
\draw (v2) [fill=white] circle (\vr);
\draw (x1) [fill=white] circle (\vr);
\draw (x2) [fill=white] circle (\vr);
\draw (x3) [fill=black] circle (\vr);
\draw (x4) [fill=white] circle (\vr);
\draw (x5) [fill=white] circle (\vr);
\draw (x6) [fill=black] circle (\vr);
\draw (y1) [fill=black] circle (\vr);
\draw (y2) [fill=white] circle (\vr);
\draw (y3) [fill=white] circle (\vr);
\draw (y4) [fill=black] circle (\vr);
\draw (y5) [fill=white] circle (\vr);
\draw (y6) [fill=white] circle (\vr);
\draw[anchor = south] (v1) node {{\small $v_1$}};
\draw[anchor = south] (v2) node {{\small $v_2$}};
\draw[anchor = north] (x1) node {{\small $x_1$}};
\draw[anchor = north] (x2) node {{\small $x_2$}};
\draw[anchor = north] (x3) node {{\small $x_3$}};
\draw[anchor = north] (x4) node {{\small $x_4$}};
\draw[anchor = north] (x5) node {{\small $x_5$}};
\draw[anchor = north] (x6) node {{\small $x_6$}};
\draw[anchor = south] (y1p) node {{\small $y_1$}};
\draw[anchor = south] (y2) node {{\small $y_2$}};
\draw[anchor = south] (y3) node {{\small $y_3$}};
\draw[anchor = south] (y4p) node {{\small $y_4$}};
\draw[anchor = south] (y5) node {{\small $y_5$}};
\draw[anchor = south] (y6) node {{\small $y_6$}};
\end{tikzpicture}
\vskip -0.25 cm
\caption{A spanning subgraph $F$ of $G$ in the proof of Claim~\ref{c:claim1}}
\label{fig:claim1}
\end{center}
\end{figure}

\begin{claim}
\label{c:claim1}
$H \ne 2C_6$.
\end{claim}
\proof Suppose, to the contrary, that $H = 2C_6$. Renaming vertices in $X$ and $Y$ if necessary, we may assume that $Q_1 \colon x_1y_1x_2y_2x_3y_3x_1$ and $Q_2 \colon x_4y_4x_5y_5x_6y_6x_4$ are the two $6$-cycles in $H$, and so $H = Q_1 \cup Q_2$. Renaming vertices if necessary, we may assume that $v_1y_1$ is an edge of $G$. Since $v_1$ is adjacent to at most two vertices from the cycle $Q_2$, we may assume, renaming vertices of $Q_2$ if necessary, that $v_2y_4$ is an edge of $G$. Thus the graph $F$ shown in Figure~\ref{fig:claim1} is a spanning subgraph of $G$. In this case, the set $S = \{y_1,x_3,y_4,x_6\}$ is a dominating set of $F$ (where the vertices in $S$ are shaded in Figure~\ref{fig:claim1}), and so $\gamma \le \gamma(F) = 4$, a contradiction.~\smallqed

\begin{figure}[ht!]
\begin{center}
\begin{tikzpicture}[scale=1,style=thick,x=1cm,y=1cm]
\def\vr{2.5pt} 
\path (0,0) coordinate (x1);
\path (0,1) coordinate (y1);
\path (-0.1,1) coordinate (y1p);
\path (1,0) coordinate (x2);
\path (1,1) coordinate (y2);
\path (2,0) coordinate (x3);
\path (2,1) coordinate (y3);
\path (1,2) coordinate (v1);
\path (3,0) coordinate (x4);
\path (3,1) coordinate (y4);
\path (2.9,1) coordinate (y4p);
\path (4,0) coordinate (x5);
\path (4,1) coordinate (y5);
\path (5,0) coordinate (x6);
\path (5,1) coordinate (y6);
\path (4,2) coordinate (v2);
\draw (x1)--(y1)--(x2)--(y2)--(x3)--(y3)--(x4)--(y4)--(x5)--(y5)--(x6)--(y6);
\draw (v1)--(y1);
\draw (v2)--(y4);
\draw (x1) to[out=-45,in=-45, distance=2.75cm] (y6);
\draw (v1) [fill=white] circle (\vr);
\draw (v2) [fill=white] circle (\vr);
\draw (x1) [fill=white] circle (\vr);
\draw (x2) [fill=white] circle (\vr);
\draw (x3) [fill=black] circle (\vr);
\draw (x4) [fill=white] circle (\vr);
\draw (x5) [fill=white] circle (\vr);
\draw (x6) [fill=black] circle (\vr);
\draw (y1) [fill=black] circle (\vr);
\draw (y2) [fill=white] circle (\vr);
\draw (y3) [fill=white] circle (\vr);
\draw (y4) [fill=black] circle (\vr);
\draw (y5) [fill=white] circle (\vr);
\draw (y6) [fill=white] circle (\vr);
\draw[anchor = south] (v1) node {{\small $v_1$}};
\draw[anchor = south] (v2) node {{\small $v_2$}};
\draw[anchor = north] (x1) node {{\small $x_1$}};
\draw[anchor = north] (x2) node {{\small $x_2$}};
\draw[anchor = north] (x3) node {{\small $x_3$}};
\draw[anchor = north] (x4) node {{\small $x_4$}};
\draw[anchor = north] (x5) node {{\small $x_5$}};
\draw[anchor = north] (x6) node {{\small $x_6$}};
\draw[anchor = south] (y1p) node {{\small $y_1$}};
\draw[anchor = south] (y2) node {{\small $y_2$}};
\draw[anchor = south] (y3) node {{\small $y_3$}};
\draw[anchor = south] (y4p) node {{\small $y_4$}};
\draw[anchor = south] (y5) node {{\small $y_5$}};
\draw[anchor = south] (y6) node {{\small $y_6$}};
\end{tikzpicture}
\vskip -1 cm
\caption{A spanning subgraph $F$ of $G$ in the proof of Claim~\ref{c:claim2}}
\label{fig:claim2}
\end{center}
\end{figure}
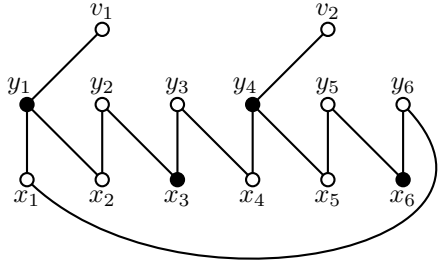

\begin{claim}
\label{c:claim2}
$H \ne C_{12}$.
\end{claim}
\proof Suppose, to the contrary, that $H = C_{12}$. Renaming vertices in $X$ and $Y$ if necessary, we may assume that $H$ is the cycle $x_1y_1x_2y_2 \ldots x_6y_6x_1$. The vertex $v_1$ has three edges to $Y$, implying that $v_1$ has exactly one edge to at least one of the three sets $\{y_1,y_4\}$, $\{y_2,y_5\}$ and $\{y_3,y_6\}$. Renaming vertices if necessary, we may assume that $v_1$ has exactly one edge to the set $\{y_1,y_4\}$. Further, we may assume that $v_1y_1$ is an edge of $G$, and so $v_1y_4$ is not an edge of $G$. Since every vertex in $Y$ is adjacent to exactly one of $v_1$ and $v_2$, this implies that $v_2y_4$ is an edge. Thus the graph $F$ shown in Figure~\ref{fig:claim2} is a spanning subgraph of $G$. In this case, the set $S = \{y_1,x_3,y_4,x_6\}$ is a dominating set of $F$ (see Figure~\ref{fig:claim2}), and so $\gamma \le \gamma(F) \le 4$, a contradiction.~\smallqed

\begin{claim}
\label{c:claim3}
If $H = 3C_{4}$, then $\pd_{\alpha}(G) \le 4$.
\end{claim}
\proof Suppose that $H = 3C_{4}$. Renaming vertices in $X$ and $Y$ if necessary, we may assume that $Q_1 \colon x_1y_1x_2y_2x_1$, $Q_2 \colon x_3y_3x_4y_4x_3$ and $Q_3 \colon x_5y_5x_6y_6x_5$ are the three $4$-cycles in $H$, and so $H = Q_1 \cup Q_2 \cup Q_2$.

\begin{figure}[ht!]
\begin{center}
\begin{tikzpicture}[scale=1,style=thick,x=1cm,y=1cm]
\def\vr{2.5pt} 
\path (0,0) coordinate (x1);
\path (0,1) coordinate (y1);
\path (-0.1,1) coordinate (y1p);
\path (1,0) coordinate (x2);
\path (1,1) coordinate (y2);
\path (0.85,1) coordinate (y2p);
\path (2,0) coordinate (x3);
\path (2,1) coordinate (y3);
\path (2.15,1) coordinate (y3p);
\path (3,0) coordinate (x4);
\path (3,1) coordinate (y4);
\path (2.9,1) coordinate (y4p);
\path (4,0) coordinate (x5);
\path (4,1) coordinate (y5);
\path (4.15,1) coordinate (y5p);
\path (5,0) coordinate (x6);
\path (5,1) coordinate (y6);
\path (1.5,2) coordinate (v1);
\path (3.5,2) coordinate (v2);
%
\draw (x1)--(y1)--(x2)--(y2)--(x1);
\draw (x3)--(y3)--(x4)--(y4)--(x3);
\draw (x5)--(y5)--(x6)--(y6)--(x5);
\draw (v1)--(y1);
\draw (v1)--(y3);
\draw (v1)--(y5);
\draw (v2)--(y2);
\draw (v2)--(y4);
\draw (v2)--(y6);
\draw (v1) [fill=black] circle (\vr);
\draw (v2) [fill=white] circle (\vr);
\draw (x1) [fill=white] circle (\vr);
\draw (x2) [fill=white] circle (\vr);
\draw (x3) [fill=white] circle (\vr);
\draw (x4) [fill=white] circle (\vr);
\draw (x5) [fill=white] circle (\vr);
\draw (x6) [fill=white] circle (\vr);
\draw (y1) [fill=white] circle (\vr);
\draw (y2) [fill=black] circle (\vr);
\draw (y3) [fill=white] circle (\vr);
\draw (y4) [fill=black] circle (\vr);
\draw (y5) [fill=white] circle (\vr);
\draw (y6) [fill=black] circle (\vr);
\draw[anchor = south] (v1) node {{\small $v_1$}};
\draw[anchor = south] (v2) node {{\small $v_2$}};
\draw[anchor = north] (x1) node {{\small $x_1$}};
\draw[anchor = north] (x2) node {{\small $x_2$}};
\draw[anchor = north] (x3) node {{\small $x_3$}};
\draw[anchor = north] (x4) node {{\small $x_4$}};
\draw[anchor = north] (x5) node {{\small $x_5$}};
\draw[anchor = north] (x6) node {{\small $x_6$}};
\draw[anchor = south] (y1p) node {{\small $y_1$}};
\draw[anchor = south] (y2p) node {{\small $y_2$}};
\draw[anchor = south] (y3p) node {{\small $y_3$}};
\draw[anchor = south] (y4p) node {{\small $y_4$}};
\draw[anchor = south] (y5p) node {{\small $y_5$}};
\draw[anchor = south] (y6) node {{\small $y_6$}};
\draw (2.5,-0.75) node {{\small (a)}};
\path (7,0) coordinate (x1);
\path (7,1) coordinate (y1);
\path (7.0,1) coordinate (y1p);
\path (8,0) coordinate (x2);
\path (8,1) coordinate (y2);
\path (7.85,1) coordinate (y2p);
\path (9,0) coordinate (x3);
\path (9,1) coordinate (y3);
\path (9.15,1) coordinate (y3p);
\path (10,0) coordinate (x4);
\path (10,1) coordinate (y4);
\path (9.9,1) coordinate (y4p);
\path (11,0) coordinate (x5);
\path (11,1) coordinate (y5);
\path (11.15,1) coordinate (y5p);
\path (12,0) coordinate (x6);
\path (12,1) coordinate (y6);
\path (8.5,2) coordinate (v1);
\path (10.5,2) coordinate (v2);
%
\draw (x1)--(y1)--(x2)--(y2)--(x1);
\draw (x3)--(y3)--(x4)--(y4)--(x3);
\draw (x5)--(y5)--(x6)--(y6)--(x5);
\draw (v1)--(y1);
\draw (v1)--(y2);
\draw (v1)--(y3);
\draw (v2)--(y4);
\draw (v2)--(y5);
\draw (v2)--(y6);
\draw (x1)--(x2);
\draw (v1) [fill=white] circle (\vr);
\draw (v2) [fill=black] circle (\vr);
\draw (x1) [fill=white] circle (\vr);
\draw (x2) [fill=black] circle (\vr);
\draw (x3) [fill=white] circle (\vr);
\draw (x4) [fill=white] circle (\vr);
\draw (x5) [fill=white] circle (\vr);
\draw (x6) [fill=white] circle (\vr);
\draw (y1) [fill=white] circle (\vr);
\draw (y2) [fill=white] circle (\vr);
\draw (y3) [fill=black] circle (\vr);
\draw (y4) [fill=white] circle (\vr);
\draw (y5) [fill=black] circle (\vr);
\draw (y6) [fill=white] circle (\vr);
\draw[anchor = south] (v1) node {{\small $v_1$}};
\draw[anchor = south] (v2) node {{\small $v_2$}};
\draw[anchor = north] (x1) node {{\small $x_1$}};
\draw[anchor = north] (x2) node {{\small $x_2$}};
\draw[anchor = north] (x3) node {{\small $x_3$}};
\draw[anchor = north] (x4) node {{\small $x_4$}};
\draw[anchor = north] (x5) node {{\small $x_5$}};
\draw[anchor = north] (x6) node {{\small $x_6$}};
\draw[anchor = south] (y1p) node {{\small $y_1$}};
\draw[anchor = south] (y2p) node {{\small $y_2$}};
\draw[anchor = south] (y3p) node {{\small $y_3$}};
\draw[anchor = south] (y4p) node {{\small $y_4$}};
\draw[anchor = south] (y5p) node {{\small $y_5$}};
\draw[anchor = south] (y6) node {{\small $y_6$}};
\draw (9.5,-0.75) node {{\small (b)}};
\end{tikzpicture}

\vskip 0.2 cm

\begin{tikzpicture}[scale=1,style=thick,x=1cm,y=1cm]
\def\vr{2.5pt} 

\path (0,0) coordinate (x1);
\path (0,1) coordinate (y1);
\path (-0.1,1) coordinate (y1p);
\path (1,0) coordinate (x2);
\path (1,1) coordinate (y2);
\path (0.85,1) coordinate (y2p);
\path (2,0) coordinate (x3);
\path (2,1) coordinate (y3);
\path (2.15,1) coordinate (y3p);
\path (3,0) coordinate (x4);
\path (3,1) coordinate (y4);
\path (2.9,1) coordinate (y4p);
\path (4,0) coordinate (x5);
\path (4,1) coordinate (y5);
\path (4.15,1) coordinate (y5p);
\path (5,0) coordinate (x6);
\path (5,1) coordinate (y6);
\path (1.5,2) coordinate (v1);
\path (3.5,2) coordinate (v2);
%
\draw (x1)--(y1)--(x2)--(y2)--(x1);
\draw (x3)--(y3)--(x4)--(y4)--(x3);
\draw (x5)--(y5)--(x6)--(y6)--(x5);
\draw (v1)--(y1);
\draw (v1)--(y2);
\draw (v1)--(y3);
\draw (v2)--(y4);
\draw (v2)--(y5);
\draw (v2)--(y6);
\draw (x1) to[out=-45,in=225, distance=1.35cm] (x6);
\draw (x2) to[out=-45,in=225, distance=0.85cm] (x5);
\draw (x3)--(x4);
\draw (v1) [fill=white] circle (\vr);
\draw (v2) [fill=white] circle (\vr);
\draw (x1) [fill=white] circle (\vr);
\draw (x2) [fill=black] circle (\vr);
\draw (x3) [fill=white] circle (\vr);
\draw (x4) [fill=white] circle (\vr);
\draw (x5) [fill=white] circle (\vr);
\draw (x6) [fill=black] circle (\vr);
\draw (y1) [fill=white] circle (\vr);
\draw (y2) [fill=white] circle (\vr);
\draw (y3) [fill=black] circle (\vr);
\draw (y4) [fill=black] circle (\vr);
\draw (y5) [fill=white] circle (\vr);
\draw (y6) [fill=white] circle (\vr);
\draw[anchor = south] (v1) node {{\small $v_1$}};
\draw[anchor = south] (v2) node {{\small $v_2$}};
\draw[anchor = north] (x1) node {{\small $x_1$}};
\draw[anchor = north] (x2) node {{\small $x_2$}};
\draw[anchor = north] (x3) node {{\small $x_3$}};
\draw[anchor = north] (x4) node {{\small $x_4$}};
\draw[anchor = north] (x5) node {{\small $x_5$}};
\draw[anchor = north] (x6) node {{\small $x_6$}};
\draw[anchor = south] (y1p) node {{\small $y_1$}};
\draw[anchor = south] (y2p) node {{\small $y_2$}};
\draw[anchor = south] (y3p) node {{\small $y_3$}};
\draw[anchor = south] (y4p) node {{\small $y_4$}};
\draw[anchor = south] (y5p) node {{\small $y_5$}};
\draw[anchor = south] (y6) node {{\small $y_6$}};
\draw (2.5,-1.25) node {{\small (c)}};
\path (7,0) coordinate (x1);
\path (7,1) coordinate (y1);
\path (6.9,1) coordinate (y1p);
\path (8,0) coordinate (x2);
\node[rectangle,draw] (y2) at (8,1) {};
\path (7.9,1.05) coordinate (y2p);
\path (9,0) coordinate (x3);
\path (9,1) coordinate (y3);
\path (9.15,1) coordinate (y3p);
\path (10,0) coordinate (x4);
\path (10,1) coordinate (y4);
\path (9.9,1) coordinate (y4p);
\path (11,0) coordinate (x5);
\path (11,1) coordinate (y5);
\path (11.15,1) coordinate (y5p);
\path (12,0) coordinate (x6);
\path (12,1) coordinate (y6);
\path (8.5,2) coordinate (v1);
\path (10.5,2) coordinate (v2);
%
\draw (x1)--(y1)--(x2)--(y2)--(x1);
\draw (x3)--(y3)--(x4)--(y4)--(x3);
\draw (x5)--(y5)--(x6)--(y6)--(x5);
\draw (v1)--(y1);
\draw (v1)--(y2);
\draw (v1)--(y3);
\draw (v2)--(y4);
\draw (v2)--(y5);
\draw (v2)--(y6);
\draw (x1) to[out=-45,in=225, distance=1.35cm] (x6);
\draw (x2)--(x3);
\draw (x4)--(x5);
\draw (v1) [fill=white] circle (\vr);
\draw (v2) [fill=black] circle (\vr);
\draw (x1) [fill=white] circle (\vr);
\draw (x2) [fill=white] circle (\vr);
\draw (x3) [fill=white] circle (\vr);
\draw (x4) [fill=white] circle (\vr);
\draw (x5) [fill=white] circle (\vr);
\draw (x6) [fill=white] circle (\vr);
\draw (y1) [fill=black] circle (\vr);
\draw (y3) [fill=black] circle (\vr);
\draw (y4) [fill=white] circle (\vr);
\draw (y5) [fill=black] circle (\vr);
\draw (y6) [fill=white] circle (\vr);
\draw[anchor = south] (v1) node {{\small $v_1$}};
\draw[anchor = south] (v2) node {{\small $v_2$}};
\draw[anchor = north] (x1) node {{\small $x_1$}};
\draw[anchor = north] (x2) node {{\small $x_2$}};
\draw[anchor = north] (x3) node {{\small $x_3$}};
\draw[anchor = north] (x4) node {{\small $x_4$}};
\draw[anchor = north] (x5) node {{\small $x_5$}};
\draw[anchor = north] (x6) node {{\small $x_6$}};
\draw[anchor = south] (y1p) node {{\small $y_1$}};
\draw[anchor = south] (y2p) node {{\small $y_2$}};
\draw[anchor = south] (y3p) node {{\small $y_3$}};
\draw[anchor = south] (y4p) node {{\small $y_4$}};
\draw[anchor = south] (y5p) node {{\small $y_5$}};
\draw[anchor = south] (y6) node {{\small $y_6$}};
\draw (9.5,-1.25) node {{\small (d) $G_{14.1}$}};
\end{tikzpicture}
\vskip -0.25 cm
\caption{Spanning subgraphs $F$ of $G$ in the proof of Claim~\ref{c:claim3}}
\label{fig:claim3}
\end{center}
\end{figure}
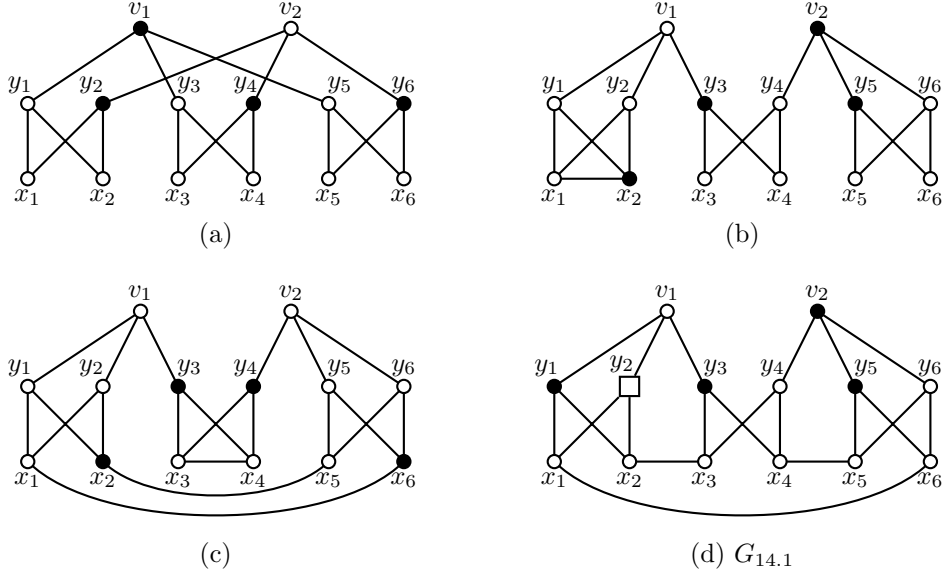

Suppose that $v_1$ is adjacent in $G$ to a vertex from each of the three $4$-cycles of $H$. Renaming vertices if necessary, we may assume that $N_G(v_1) = \{y_1,y_3,y_5\}$. Since every vertex in $Y$ is adjacent to exactly one of $v_1$ and $v_2$, this implies that $N_G(v_2) = \{y_2,y_4,y_6\}$. Thus the graph $F$ shown in Figure~\ref{fig:claim3}(a) is a spanning subgraph of $G$. In this case, the set $S = \{v_1,y_2,y_4,y_6\}$ is a dominating set of $F$ (see Figure~\ref{fig:claim3}(a)), and so $\gamma \le \gamma(F) = 4$, a contradiction.

Hence, neither $v_1$ nor $v_2$ is adjacent in $G$ to a vertex from each of the three $4$-cycles of $H$. Renaming vertices if necessary, we may assume that $N_G(v_1) = \{y_1,y_2,y_3\}$ and $N_G(v_2) = \{y_4,y_5,y_6\}$. By our earlier observations, $G[X] = 3P_2$. If $x_1x_2$ is an edge, then the graph $F$ shown in Figure~\ref{fig:claim3}(b) is a spanning subgraph of $G$. In this case, the set $S = \{v_2,x_2,y_3,y_5\}$ is a dominating set of $F$ (see Figure~\ref{fig:claim3}(b)), and so $\gamma \le \gamma(F) = 4$, a contradiction. Hence, $x_1x_2 \notin E(G)$. By symmetry, $x_5x_6 \notin E(G)$.

Suppose that $x_3x_4$ is an edge. Renaming vertices in necessary, we may assume in this case that $x_1x_6$ and $x_2x_5$ are edges. Thus the graph $G$ is determined and is shown in Figure~\ref{fig:claim3}(c). In this case, the set $S = \{x_2,x_6,y_3,y_4\}$ is a dominating set of $G$ (see Figure~\ref{fig:claim3}(c)), and so $\gamma \le 4$, a contradiction. Hence, $x_3x_4 \notin E(G)$.

The graph $G$ is therefore determined. Renaming vertices if necessary, we may assume that $G = G_{14.1}$, where $G_{14.1}$  is the graph shown in Figure~\ref{fig:claim3}(d). We note that $\gamma = 5$. In this case, the set $S = \{y_1,y_3,y_5,v_2\}$ satisfies $\dom_G(S) = 13$ (the vertex $y_2$ represented by the square in Figure~\ref{fig:claim3}(d) is the only vertex not dominated by $S$) and $|S| = 4$, implying that $\pd_{\alpha}(G) \le |S| \le 4$. This completes the proof of Claim~\ref{c:claim3}.~\smallqed

\begin{claim}
\label{c:claim4}
If $H = C_{4} \cup C_8$, then $\pd_{\alpha}(G) \le 4$.
\end{claim}
\proof Suppose that $H = C_{4} \cup C_8$. Renaming vertices in $X$ and $Y$ if necessary, we may assume that $Q_1 \colon x_1y_1x_2y_2x_1$ is the $4$-cycle in $H$ and $Q_2 \colon x_3y_3x_4y_4x_5y_5x_6y_6x_3$ is the $6$-cycle in $H$.

\begin{subclaim}
\label{c:claim4.1}
Both $v_1$ and $v_2$ are adjacent to exactly one vertex from the cycle $Q_1$.
\end{subclaim}
\proof Suppose, to the contrary, that $v_1$ or $v_2$, say $v_1$, is adjacent in $G$ to two vertices in the cycle $Q_1$. Renaming vertices if necessary, we may assume that $N_G(v_1) = \{y_1,y_2,y_3\}$, and so $N_G(v_2) = \{y_4,y_5,y_6\}$. Recall that  $G[X] = 3P_2$.

\begin{figure}[ht!]
\begin{center}
\begin{tikzpicture}[scale=1,style=thick,x=1cm,y=1cm]
\def\vr{2.5pt} 
\path (0,0) coordinate (x1);
\path (0,1) coordinate (y1);
\path (-0.1,1) coordinate (y1p);
\path (1,0) coordinate (x2);
\path (1,1) coordinate (y2);
\path (0.85,1) coordinate (y2p);
\path (2,0) coordinate (x3);
\path (2,1) coordinate (y3);
\path (2.15,1) coordinate (y3p);
\path (3,0) coordinate (x4);
\path (3,1) coordinate (y4);
\path (2.9,1) coordinate (y4p);
\path (4,0) coordinate (x5);
\path (4,1) coordinate (y5);
\path (4.15,1) coordinate (y5p);
\path (5,0) coordinate (x6);
\path (5,1) coordinate (y6);
\path (1.5,2) coordinate (v1);
\path (3.5,2) coordinate (v2);
%
\draw (x1)--(y1)--(x2)--(y2)--(x1);
\draw (x3)--(y3)--(x4)--(y4)--(x5)--(y5)--(x6)--(y6)--(x3);
\draw (v1)--(y1);
\draw (v1)--(y2);
\draw (v1)--(y3);
\draw (v2)--(y4);
\draw (v2)--(y5);
\draw (v2)--(y6);
\draw (x1)--(x2);
\draw (v1) [fill=white] circle (\vr);
\draw (v2) [fill=white] circle (\vr);
\draw (x1) [fill=black] circle (\vr);
\draw (x2) [fill=white] circle (\vr);
\draw (x3) [fill=white] circle (\vr);
\draw (x4) [fill=white] circle (\vr);
\draw (x5) [fill=white] circle (\vr);
\draw (x6) [fill=black] circle (\vr);
\draw (y1) [fill=white] circle (\vr);
\draw (y2) [fill=white] circle (\vr);
\draw (y3) [fill=black] circle (\vr);
\draw (y4) [fill=black] circle (\vr);
\draw (y5) [fill=white] circle (\vr);
\draw (y6) [fill=white] circle (\vr);
\draw[anchor = south] (v1) node {{\small $v_1$}};
\draw[anchor = south] (v2) node {{\small $v_2$}};
\draw[anchor = north] (x1) node {{\small $x_1$}};
\draw[anchor = north] (x2) node {{\small $x_2$}};
\draw[anchor = north] (x3) node {{\small $x_3$}};
\draw[anchor = north] (x4) node {{\small $x_4$}};
\draw[anchor = north] (x5) node {{\small $x_5$}};
\draw[anchor = north] (x6) node {{\small $x_6$}};
\draw[anchor = south] (y1p) node {{\small $y_1$}};
\draw[anchor = south] (y2p) node {{\small $y_2$}};
\draw[anchor = south] (y3p) node {{\small $y_3$}};
\draw[anchor = south] (y4p) node {{\small $y_4$}};
\draw[anchor = south] (y5p) node {{\small $y_5$}};
\draw[anchor = south] (y6) node {{\small $y_6$}};
\draw (2.5,-0.75) node {{\small (a)}};
\path (7,0) coordinate (x1);
\path (7,1) coordinate (y1);
\path (6.9,1) coordinate (y1p);
\path (8,0) coordinate (x2);
\path (8,1) coordinate (y2);
\path (7.85,1) coordinate (y2p);
\path (9,0) coordinate (x3);
\path (9,1) coordinate (y3);
\path (9.15,1) coordinate (y3p);
\path (10,0) coordinate (x4);
\path (10.1,0) coordinate (x4p);
\path (10,1) coordinate (y4);
\path (9.9,1) coordinate (y4p);
\path (11,0) coordinate (x5);
\path (11,1) coordinate (y5);
\path (11.15,1) coordinate (y5p);
\path (12,0) coordinate (x6);
\path (12,1) coordinate (y6);
\path (8.5,2) coordinate (v1);
\path (10.5,2) coordinate (v2);
%
\draw (x1)--(y1)--(x2)--(y2)--(x1);
\draw (x3)--(y3)--(x4)--(y4)--(x5)--(y5)--(x6)--(y6)--(x3);
\draw (v1)--(y1);
\draw (v1)--(y2);
\draw (v1)--(y3);
\draw (v2)--(y4);
\draw (v2)--(y5);
\draw (v2)--(y6);
\draw (x2) to[out=-45,in=225, distance=0.95cm] (x6);
\draw (v1) [fill=white] circle (\vr);
\draw (v2) [fill=white] circle (\vr);
\draw (x1) [fill=black] circle (\vr);
\draw (x2) [fill=white] circle (\vr);
\draw (x3) [fill=white] circle (\vr);
\draw (x4) [fill=white] circle (\vr);
\draw (x5) [fill=white] circle (\vr);
\draw (x6) [fill=black] circle (\vr);
\draw (y1) [fill=white] circle (\vr);
\draw (y2) [fill=white] circle (\vr);
\draw (y3) [fill=black] circle (\vr);
\draw (y4) [fill=black] circle (\vr);
\draw (y5) [fill=white] circle (\vr);
\draw (y6) [fill=white] circle (\vr);
\draw[anchor = south] (v1) node {{\small $v_1$}};
\draw[anchor = south] (v2) node {{\small $v_2$}};
\draw[anchor = north] (x1) node {{\small $x_1$}};
\draw[anchor = north] (x2) node {{\small $x_2$}};
\draw[anchor = north] (x3) node {{\small $x_3$}};
\draw[anchor = north] (x4) node {{\small $x_4$}};
\draw[anchor = north] (x5) node {{\small $x_5$}};
\draw[anchor = north] (x6) node {{\small $x_6$}};
\draw[anchor = south] (y1p) node {{\small $y_1$}};
\draw[anchor = south] (y2p) node {{\small $y_2$}};
\draw[anchor = south] (y3p) node {{\small $y_3$}};
\draw[anchor = south] (y4p) node {{\small $y_4$}};
\draw[anchor = south] (y5p) node {{\small $y_5$}};
\draw[anchor = south] (y6) node {{\small $y_6$}};
\draw (9.5,-0.75) node {{\small (b)}};
\end{tikzpicture}

\vskip 0.2 cm

\begin{tikzpicture}[scale=1,style=thick,x=1cm,y=1cm]
\def\vr{2.5pt} 
\path (0,0) coordinate (x1);
\path (0,1) coordinate (y1);
\path (-0.1,1) coordinate (y1p);
\path (1,0) coordinate (x2);
\path (1,1) coordinate (y2);
\path (0.85,1) coordinate (y2p);
\path (2,0) coordinate (x3);
\path (2,1) coordinate (y3);
\path (2.15,1) coordinate (y3p);
\path (3,0) coordinate (x4);
\path (3,1) coordinate (y4);
\path (2.9,1) coordinate (y4p);
\path (4,0) coordinate (x5);
\path (4.05,0) coordinate (x5p);
\path (4,1) coordinate (y5);
\path (4.15,1) coordinate (y5p);
\path (5,0) coordinate (x6);
\path (5,1) coordinate (y6);
\path (1.5,2) coordinate (v1);
\path (3.5,2) coordinate (v2);
%
\draw (x1)--(y1)--(x2)--(y2)--(x1);
\draw (x3)--(y3)--(x4)--(y4)--(x5)--(y5)--(x6)--(y6)--(x3);
\draw (v1)--(y1);
\draw (v1)--(y2);
\draw (v1)--(y3);
\draw (v2)--(y4);
\draw (v2)--(y5);
\draw (v2)--(y6);
\draw (x2) to[out=-45,in=225, distance=0.9cm] (x5);
\draw (v1) [fill=white] circle (\vr);
\draw (v2) [fill=white] circle (\vr);
\draw (x1) [fill=black] circle (\vr);
\draw (x2) [fill=white] circle (\vr);
\draw (x3) [fill=white] circle (\vr);
\draw (x4) [fill=white] circle (\vr);
\draw (x5) [fill=black] circle (\vr);
\draw (x6) [fill=white] circle (\vr);
\draw (y1) [fill=white] circle (\vr);
\draw (y2) [fill=white] circle (\vr);
\draw (y3) [fill=black] circle (\vr);
\draw (y4) [fill=white] circle (\vr);
\draw (y5) [fill=white] circle (\vr);
\draw (y6) [fill=black] circle (\vr);
\draw[anchor = south] (v1) node {{\small $v_1$}};
\draw[anchor = south] (v2) node {{\small $v_2$}};
\draw[anchor = north] (x1) node {{\small $x_1$}};
\draw[anchor = north] (x2) node {{\small $x_2$}};
\draw[anchor = north] (x3) node {{\small $x_3$}};
\draw[anchor = north] (x4) node {{\small $x_4$}};
\draw[anchor = north] (x5p) node {{\small $x_5$}};
\draw[anchor = north] (x6) node {{\small $x_6$}};
\draw[anchor = south] (y1p) node {{\small $y_1$}};
\draw[anchor = south] (y2p) node {{\small $y_2$}};
\draw[anchor = south] (y3p) node {{\small $y_3$}};
\draw[anchor = south] (y4p) node {{\small $y_4$}};
\draw[anchor = south] (y5p) node {{\small $y_5$}};
\draw[anchor = south] (y6) node {{\small $y_6$}};
\draw (2.5,-0.75) node {{\small (c)}};
\path (7,0) coordinate (x1);
\path (7,1) coordinate (y1);
\path (6.9,1) coordinate (y1p);
\path (8,0) coordinate (x2);
\path (8,1) coordinate (y2);
\path (7.85,1) coordinate (y2p);
\path (9,0) coordinate (x3);
\path (9,1) coordinate (y3);
\path (9.15,1) coordinate (y3p);
\path (10,0) coordinate (x4);
\path (10.05,0) coordinate (x4p);
\path (10,1) coordinate (y4);
\path (9.9,1) coordinate (y4p);
\path (11,0) coordinate (x5);
\path (11,1) coordinate (y5);
\path (11.15,1) coordinate (y5p);
\path (12,0) coordinate (x6);
\path (12,1) coordinate (y6);
\path (8.5,2) coordinate (v1);
\path (10.5,2) coordinate (v2);
%
\draw (x1)--(y1)--(x2)--(y2)--(x1);
\draw (x3)--(y3)--(x4)--(y4)--(x5)--(y5)--(x6)--(y6)--(x3);
\draw (v1)--(y1);
\draw (v1)--(y2);
\draw (v1)--(y3);
\draw (v2)--(y4);
\draw (v2)--(y5);
\draw (v2)--(y6);
\draw (x1) to[out=-45,in=225, distance=0.95cm] (x4);
\draw (x2)--(x3);
\draw (x5)--(x6);
\draw (v1) [fill=black] circle (\vr);
\draw (v2) [fill=white] circle (\vr);
\draw (x1) [fill=white] circle (\vr);
\draw (x2) [fill=white] circle (\vr);
\draw (x3) [fill=black] circle (\vr);
\draw (x4) [fill=black] circle (\vr);
\draw (x5) [fill=white] circle (\vr);
\draw (x6) [fill=white] circle (\vr);
\draw (y1) [fill=white] circle (\vr);
\draw (y2) [fill=white] circle (\vr);
\draw (y3) [fill=white] circle (\vr);
\draw (y4) [fill=white] circle (\vr);
\draw (y5) [fill=black] circle (\vr);
\draw (y6) [fill=white] circle (\vr);
\draw[anchor = south] (v1) node {{\small $v_1$}};
\draw[anchor = south] (v2) node {{\small $v_2$}};
\draw[anchor = north] (x1) node {{\small $x_1$}};
\draw[anchor = north] (x2) node {{\small $x_2$}};
\draw[anchor = north] (x3) node {{\small $x_3$}};
\draw[anchor = north] (x4p) node {{\small $x_4$}};
\draw[anchor = north] (x5) node {{\small $x_5$}};
\draw[anchor = north] (x6) node {{\small $x_6$}};
\draw[anchor = south] (y1p) node {{\small $y_1$}};
\draw[anchor = south] (y2p) node {{\small $y_2$}};
\draw[anchor = south] (y3p) node {{\small $y_3$}};
\draw[anchor = south] (y4p) node {{\small $y_4$}};
\draw[anchor = south] (y5p) node {{\small $y_5$}};
\draw[anchor = south] (y6) node {{\small $y_6$}};
\draw (9.5,-0.75) node {{\small (d)}};
\end{tikzpicture}
\vskip -0.25 cm
\caption{Spanning subgraphs $F$ of $G$ in the proof of Claim~\ref{c:claim4.1}}
\label{fig:claim4.1}
\end{center}
\end{figure}
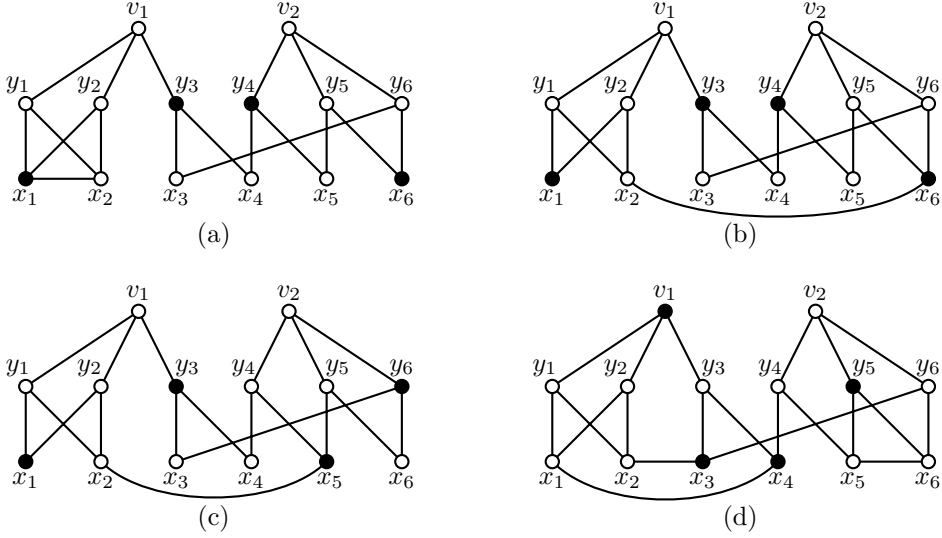

If $x_1x_2$ is an edge, then the graph $F$ shown in Figure~\ref{fig:claim4.1}(a) is a spanning subgraph of $G$. In this case, the set $S = \{x_1,x_6,y_3,y_4\}$ is a dominating set of $F$ (see Figure~\ref{fig:claim4.1}(a)), and so $\gamma \le \gamma(F) = 4$, a contradiction. Hence, $x_1x_2 \notin E(G)$.
If $x_2x_6$ is an edge, then the graph $F$ shown in Figure~\ref{fig:claim4.1}(b) is a spanning subgraph of $G$. In this case, the set $S = \{x_1,x_6,y_3,y_4\}$ is a dominating set of $F$ (see Figure~\ref{fig:claim4.1}(b)), and so $\gamma \le \gamma(F) = 4$, a contradiction. Hence, $x_2x_6 \notin E(G)$. By symmetry, $x_1x_6 \notin E(G)$.
If $x_2x_5$ is an edge, then the graph $F$ shown in Figure~\ref{fig:claim4.1}(c) is a spanning subgraph of $G$. In this case, the set $S = \{x_1,x_5,y_3,y_6\}$ is a dominating set of $F$ (see Figure~\ref{fig:claim4.1}(c)), and so $\gamma \le \gamma(F) = 4$, a contradiction. Hence, $x_2x_5 \notin E(G)$. By symmetry, $x_1x_5 \notin E(G)$.

Renaming $x_1$ and $x_2$ if necessary, we may assume that $x_1x_4$ and $x_2x_3$ are edges. The remaining edge in $G[X]$ is therefore the edge $x_5x_6$. Thus, the graph $G$ is determined and is shown in Figure~\ref{fig:claim4.1}(d). In this case, the set $S = \{v_1,x_3,x_4,y_5\}$ is a dominating set of $G$ (see Figure~\ref{fig:claim4.1}(d)), and so $\gamma \le 4$, a contradiction. This completes the proof of Claim~\ref{c:claim4.1}.~\smallqed

\medskip
By Claim~\ref{c:claim4.1}, both $v_1$ and $v_2$ are adjacent to exactly one vertex from the cycle $Q_1$. Renaming $y_1$ and $y_2$ if necessary, we may assume that $v_1y_1$ and $v_2y_2$ are edges.

\begin{figure}[ht!]
\begin{center}
\begin{tikzpicture}[scale=1,style=thick,x=1cm,y=1cm]
\def\vr{2.5pt} 
\path (0,0) coordinate (x1);
\path (0,1) coordinate (y1);
\path (-0.1,1) coordinate (y1p);
\path (1,0) coordinate (x2);
\path (1,1) coordinate (y2);
\path (0.85,1) coordinate (y2p);
\path (2,0) coordinate (x3);
\path (2,1) coordinate (y3);
\path (2.15,1) coordinate (y3p);
\path (3,0) coordinate (x4);
\path (3,1) coordinate (y4);
\path (2.9,1) coordinate (y4p);
\path (4,0) coordinate (x5);
\path (4,1) coordinate (y5);
\path (4.15,1) coordinate (y5p);
\path (5,0) coordinate (x6);
\path (5,1) coordinate (y6);
\path (1.5,2) coordinate (v1);
\path (3.5,2) coordinate (v2);
%
\draw (x1)--(y1)--(x2)--(y2)--(x1);
\draw (x3)--(y3)--(x4)--(y4)--(x5)--(y5)--(x6)--(y6)--(x3);
\draw (v1)--(y1);
\draw (v1)--(y3);
\draw (v1)--(y5);
\draw (v2)--(y2);
\draw (v2)--(y4);
\draw (v2)--(y6);
\draw (v1) [fill=black] circle (\vr);
\draw (v2) [fill=white] circle (\vr);
\draw (x1) [fill=white] circle (\vr);
\draw (x2) [fill=white] circle (\vr);
\draw (x3) [fill=white] circle (\vr);
\draw (x4) [fill=white] circle (\vr);
\draw (x5) [fill=white] circle (\vr);
\draw (x6) [fill=white] circle (\vr);
\draw (y1) [fill=white] circle (\vr);
\draw (y2) [fill=black] circle (\vr);
\draw (y3) [fill=white] circle (\vr);
\draw (y4) [fill=black] circle (\vr);
\draw (y5) [fill=white] circle (\vr);
\draw (y6) [fill=black] circle (\vr);
\draw[anchor = south] (v1) node {{\small $v_1$}};
\draw[anchor = south] (v2) node {{\small $v_2$}};
\draw[anchor = north] (x1) node {{\small $x_1$}};
\draw[anchor = north] (x2) node {{\small $x_2$}};
\draw[anchor = north] (x3) node {{\small $x_3$}};
\draw[anchor = north] (x4) node {{\small $x_4$}};
\draw[anchor = north] (x5) node {{\small $x_5$}};
\draw[anchor = north] (x6) node {{\small $x_6$}};
\draw[anchor = south] (y1p) node {{\small $y_1$}};
\draw[anchor = south] (y2p) node {{\small $y_2$}};
\draw[anchor = south] (y3p) node {{\small $y_3$}};
\draw[anchor = south] (y4p) node {{\small $y_4$}};
\draw[anchor = south] (y5p) node {{\small $y_5$}};
\draw[anchor = south] (y6) node {{\small $y_6$}};
\end{tikzpicture}
\vskip -0.25 cm
\caption{A spanning subgraph $F$ of $G$ in the proof of Claim~\ref{c:claim4.2}}
\label{fig:claim4.2}
\end{center}
\end{figure}
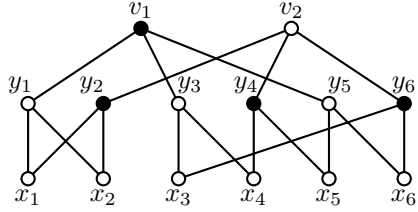

\begin{subclaim}
\label{c:claim4.2}
The vertex $v_1$ is adjacent to two vertices in $Q_2$ at distance~$2$ in $Q_2$.
\end{subclaim}
\proof Suppose, to the contrary, that $v_1$ is adjacent to two vertices in $Q_2$ at distance~$4$. Renaming vertices if necessary, we may assume that $v_1y_3$ and $v_1y_5$ are edges. Thus, $N_G(v_1) = \{y_1,y_3,y_5\}$ and $N_G(v_2) = \{y_2,y_4,y_6\}$. Thus the graph $F$ shown in Figure~\ref{fig:claim4.2} is a spanning subgraph of $G$. In this case, the set $S = \{v_1,y_2,y_4,y_6\}$ is a dominating set of $F$ (see Figure~\ref{fig:claim4.2}), and so $\gamma \le \gamma(F) = 4$, a contradiction.~\smallqed

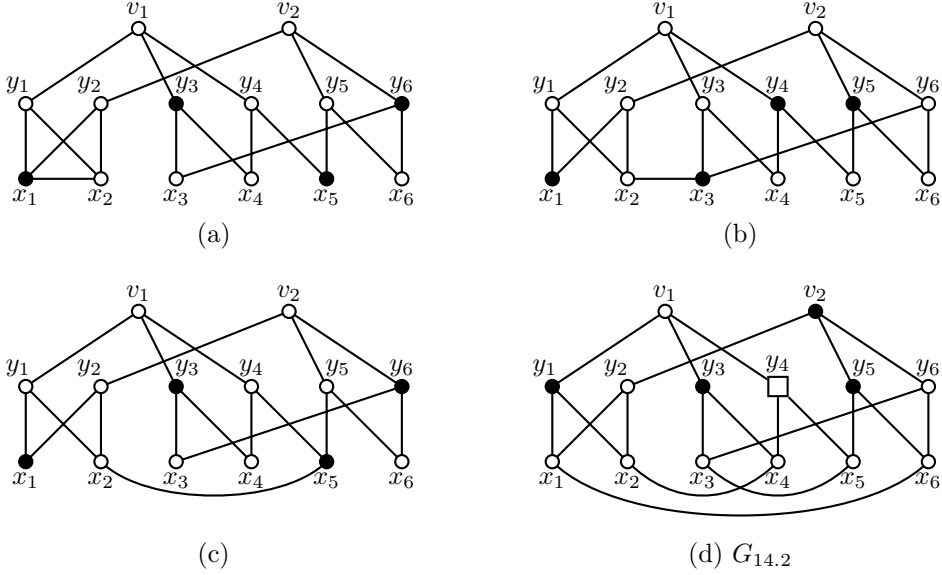
\begin{figure}[ht!]
\begin{center}
\begin{tikzpicture}[scale=1,style=thick,x=1cm,y=1cm]
\def\vr{2.5pt} 
\path (0,0) coordinate (x1);
\path (0,1) coordinate (y1);
\path (-0.1,1) coordinate (y1p);
\path (1,0) coordinate (x2);
\path (1,1) coordinate (y2);
\path (0.85,1) coordinate (y2p);
\path (2,0) coordinate (x3);
\path (2,1) coordinate (y3);
\path (2.15,1) coordinate (y3p);
\path (3,0) coordinate (x4);
\path (3,1) coordinate (y4);
\path (3.1,1) coordinate (y4p);
\path (4,0) coordinate (x5);
\path (4,1) coordinate (y5);
\path (4.15,1) coordinate (y5p);
\path (5,0) coordinate (x6);
\path (5,1) coordinate (y6);
\path (1.5,2) coordinate (v1);
\path (3.5,2) coordinate (v2);
%
\draw (x1)--(y1)--(x2)--(y2)--(x1);
\draw (x3)--(y3)--(x4)--(y4)--(x5)--(y5)--(x6)--(y6)--(x3);
\draw (v1)--(y1);
\draw (v1)--(y3);
\draw (v1)--(y4);
\draw (v2)--(y2);
\draw (v2)--(y5);
\draw (v2)--(y6);
\draw (x1)--(x2);
\draw (v1) [fill=white] circle (\vr);
\draw (v2) [fill=white] circle (\vr);
\draw (x1) [fill=black] circle (\vr);
\draw (x2) [fill=white] circle (\vr);
\draw (x3) [fill=white] circle (\vr);
\draw (x4) [fill=white] circle (\vr);
\draw (x5) [fill=black] circle (\vr);
\draw (x6) [fill=white] circle (\vr);
\draw (y1) [fill=white] circle (\vr);
\draw (y2) [fill=white] circle (\vr);
\draw (y3) [fill=black] circle (\vr);
\draw (y4) [fill=white] circle (\vr);
\draw (y5) [fill=white] circle (\vr);
\draw (y6) [fill=black] circle (\vr);
\draw[anchor = south] (v1) node {{\small $v_1$}};
\draw[anchor = south] (v2) node {{\small $v_2$}};
\draw[anchor = north] (x1) node {{\small $x_1$}};
\draw[anchor = north] (x2) node {{\small $x_2$}};
\draw[anchor = north] (x3) node {{\small $x_3$}};
\draw[anchor = north] (x4) node {{\small $x_4$}};
\draw[anchor = north] (x5) node {{\small $x_5$}};
\draw[anchor = north] (x6) node {{\small $x_6$}};
\draw[anchor = south] (y1p) node {{\small $y_1$}};
\draw[anchor = south] (y2p) node {{\small $y_2$}};
\draw[anchor = south] (y3p) node {{\small $y_3$}};
\draw[anchor = south] (y4) node {{\small $y_4$}};
\draw[anchor = south] (y5p) node {{\small $y_5$}};
\draw[anchor = south] (y6) node {{\small $y_6$}};
\draw (2.5,-0.75) node {{\small (a)}};
\path (7,0) coordinate (x1);
\path (7,1) coordinate (y1);
\path (6.9,1) coordinate (y1p);
\path (8,0) coordinate (x2);
\path (8,1) coordinate (y2);
\path (7.85,1) coordinate (y2p);
\path (9,0) coordinate (x3);
\path (9,1) coordinate (y3);
\path (9.15,1) coordinate (y3p);
\path (10,0) coordinate (x4);
\path (10,1) coordinate (y4);
\path (10.1,1) coordinate (y4p);
\path (11,0) coordinate (x5);
\path (11,1) coordinate (y5);
\path (11.15,1) coordinate (y5p);
\path (12,0) coordinate (x6);
\path (12,1) coordinate (y6);
\path (8.5,2) coordinate (v1);
\path (10.5,2) coordinate (v2);
%
\draw (x1)--(y1)--(x2)--(y2)--(x1);
\draw (x3)--(y3)--(x4)--(y4)--(x5)--(y5)--(x6)--(y6)--(x3);
\draw (v1)--(y1);
\draw (v1)--(y3);
\draw (v1)--(y4);
\draw (v2)--(y2);
\draw (v2)--(y5);
\draw (v2)--(y6);
\draw (x2)--(x3);
\draw (v1) [fill=white] circle (\vr);
\draw (v2) [fill=white] circle (\vr);
\draw (x1) [fill=black] circle (\vr);
\draw (x2) [fill=white] circle (\vr);
\draw (x3) [fill=black] circle (\vr);
\draw (x4) [fill=white] circle (\vr);
\draw (x5) [fill=white] circle (\vr);
\draw (x6) [fill=white] circle (\vr);
\draw (y1) [fill=white] circle (\vr);
\draw (y2) [fill=white] circle (\vr);
\draw (y3) [fill=white] circle (\vr);
\draw (y4) [fill=black] circle (\vr);
\draw (y5) [fill=black] circle (\vr);
\draw (y6) [fill=white] circle (\vr);
\draw[anchor = south] (v1) node {{\small $v_1$}};
\draw[anchor = south] (v2) node {{\small $v_2$}};
\draw[anchor = north] (x1) node {{\small $x_1$}};
\draw[anchor = north] (x2) node {{\small $x_2$}};
\draw[anchor = north] (x3) node {{\small $x_3$}};
\draw[anchor = north] (x4) node {{\small $x_4$}};
\draw[anchor = north] (x5) node {{\small $x_5$}};
\draw[anchor = north] (x6) node {{\small $x_6$}};
\draw[anchor = south] (y1p) node {{\small $y_1$}};
\draw[anchor = south] (y2p) node {{\small $y_2$}};
\draw[anchor = south] (y3p) node {{\small $y_3$}};
\draw[anchor = south] (y4) node {{\small $y_4$}};
\draw[anchor = south] (y5p) node {{\small $y_5$}};
\draw[anchor = south] (y6) node {{\small $y_6$}};
\draw (9.5,-0.75) node {{\small (b)}};
\end{tikzpicture}

\vskip 0.2 cm

\begin{tikzpicture}[scale=1,style=thick,x=1cm,y=1cm]
\def\vr{2.5pt} 
\path (0,0) coordinate (x1);
\path (0,1) coordinate (y1);
\path (-0.1,1) coordinate (y1p);
\path (1,0) coordinate (x2);
\path (1,1) coordinate (y2);
\path (0.85,1) coordinate (y2p);
\path (2,0) coordinate (x3);
\path (2,1) coordinate (y3);
\path (2.15,1) coordinate (y3p);
\path (3,0) coordinate (x4);
\path (3,1) coordinate (y4);
\path (3.1,1) coordinate (y4p);
\path (4,0) coordinate (x5);
\path (4,1) coordinate (y5);
\path (4.15,1) coordinate (y5p);
\path (5,0) coordinate (x6);
\path (5,1) coordinate (y6);
\path (1.5,2) coordinate (v1);
\path (3.5,2) coordinate (v2);
%
\draw (x1)--(y1)--(x2)--(y2)--(x1);
\draw (x3)--(y3)--(x4)--(y4)--(x5)--(y5)--(x6)--(y6)--(x3);
\draw (v1)--(y1);
\draw (v1)--(y3);
\draw (v1)--(y4);
\draw (v2)--(y2);
\draw (v2)--(y5);
\draw (v2)--(y6);
\draw (x2) to[out=-45,in=225, distance=0.85cm] (x5);
\draw (v1) [fill=white] circle (\vr);
\draw (v2) [fill=white] circle (\vr);
\draw (x1) [fill=black] circle (\vr);
\draw (x2) [fill=white] circle (\vr);
\draw (x3) [fill=white] circle (\vr);
\draw (x4) [fill=white] circle (\vr);
\draw (x5) [fill=black] circle (\vr);
\draw (x6) [fill=white] circle (\vr);
\draw (y1) [fill=white] circle (\vr);
\draw (y2) [fill=white] circle (\vr);
\draw (y3) [fill=black] circle (\vr);
\draw (y4) [fill=white] circle (\vr);
\draw (y5) [fill=white] circle (\vr);
\draw (y6) [fill=black] circle (\vr);
\draw[anchor = south] (v1) node {{\small $v_1$}};
\draw[anchor = south] (v2) node {{\small $v_2$}};
\draw[anchor = north] (x1) node {{\small $x_1$}};
\draw[anchor = north] (x2) node {{\small $x_2$}};
\draw[anchor = north] (x3) node {{\small $x_3$}};
\draw[anchor = north] (x4) node {{\small $x_4$}};
\draw[anchor = north] (x5) node {{\small $x_5$}};
\draw[anchor = north] (x6) node {{\small $x_6$}};
\draw[anchor = south] (y1p) node {{\small $y_1$}};
\draw[anchor = south] (y2p) node {{\small $y_2$}};
\draw[anchor = south] (y3p) node {{\small $y_3$}};
\draw[anchor = south] (y4) node {{\small $y_4$}};
\draw[anchor = south] (y5p) node {{\small $y_5$}};
\draw[anchor = south] (y6) node {{\small $y_6$}};
\draw (2.5,-1.25) node {{\small (c)}};
\path (7,0) coordinate (x1);
\path (7,1) coordinate (y1);
\path (6.9,1) coordinate (y1p);
\path (8,0) coordinate (x2);
\path (8,1) coordinate (y2);
\path (7.85,1) coordinate (y2p);
\path (9,0) coordinate (x3);
\path (9,1) coordinate (y3);
\path (9.15,1) coordinate (y3p);
\path (10,0) coordinate (x4);
\node[rectangle,draw] (y4) at (10,1) {};
\path (10,1.05) coordinate (y4p);
\path (11,0) coordinate (x5);
\path (11,1) coordinate (y5);
\path (11.15,1) coordinate (y5p);
\path (12,0) coordinate (x6);
\path (12,1) coordinate (y6);
\path (8.5,2) coordinate (v1);
\path (10.5,2) coordinate (v2);
%
\draw (x1)--(y1)--(x2)--(y2)--(x1);
\draw (x3)--(y3)--(x4)--(y4)--(x5)--(y5)--(x6)--(y6)--(x3);
\draw (v1)--(y1);
\draw (v1)--(y3);
\draw (v1)--(y4);
\draw (v2)--(y2);
\draw (v2)--(y5);
\draw (v2)--(y6);
\draw (x1) to[out=-45,in=225, distance=1.35cm] (x6);
\draw (x2) to[out=-45,in=225, distance=0.85cm] (x4);
\draw (x3) to[out=-45,in=225, distance=0.85cm] (x5);
\draw (v1) [fill=white] circle (\vr);
\draw (v2) [fill=black] circle (\vr);
\draw (x1) [fill=white] circle (\vr);
\draw (x2) [fill=white] circle (\vr);
\draw (x3) [fill=white] circle (\vr);
\draw (x4) [fill=white] circle (\vr);
\draw (x5) [fill=white] circle (\vr);
\draw (x6) [fill=white] circle (\vr);
\draw (y1) [fill=black] circle (\vr);
\draw (y2) [fill=white] circle (\vr);
\draw (y3) [fill=black] circle (\vr);
\draw (y5) [fill=black] circle (\vr);
\draw (y6) [fill=white] circle (\vr);
\draw[anchor = south] (v1) node {{\small $v_1$}};
\draw[anchor = south] (v2) node {{\small $v_2$}};
\draw[anchor = north] (x1) node {{\small $x_1$}};
\draw[anchor = north] (x2) node {{\small $x_2$}};
\draw[anchor = north] (x3) node {{\small $x_3$}};
\draw[anchor = north] (x4) node {{\small $x_4$}};
\draw[anchor = north] (x5) node {{\small $x_5$}};
\draw[anchor = north] (x6) node {{\small $x_6$}};
\draw[anchor = south] (y1p) node {{\small $y_1$}};
\draw[anchor = south] (y2p) node {{\small $y_2$}};
\draw[anchor = south] (y3p) node {{\small $y_3$}};
\draw[anchor = south] (y4p) node {{\small $y_4$}};
\draw[anchor = south] (y5p) node {{\small $y_5$}};
\draw[anchor = south] (y6) node {{\small $y_6$}};
\draw (9.5,-1.25) node {{\small (d) $G_{14.2}$}};
\end{tikzpicture}
\vskip -0.25 cm
\caption{Spanning subgraphs $F$ of $G$ in the proof of Claim~\ref{c:claim4}}
\label{fig:claim4sub}
\end{center}
\end{figure}

\medskip
By Claim~\ref{c:claim4.2}, the vertex $v_1$ is adjacent to two vertices in $Q_2$ at distance~$2$ in $Q_2$. Renaming vertices if necessary, we may assume that $N_G(v_1) = \{y_1,y_3,y_4\}$ and $N_G(v_2) = \{y_2,y_5,y_6\}$. If $x_1x_2$ is an edge, then the graph $F$ shown in Figure~\ref{fig:claim4sub}(a) is a spanning subgraph of $G$. In this case, the set $S = \{x_1,x_5,y_3,y_6\}$ is a dominating set of $F$ (see Figure~\ref{fig:claim4sub}(a)), and so $\gamma \le \gamma(F) = 4$, a contradiction. Hence, $x_1x_2 \notin E(G)$.

If $x_2x_3$ is an edge, then the graph $F$ shown in Figure~\ref{fig:claim4sub}(b) is a spanning subgraph of $G$. In this case, the set $S = \{x_1,x_3,y_4,y_5\}$ is a dominating set of $F$ (see Figure~\ref{fig:claim4sub}(b)), and so $\gamma \le \gamma(F) = 4$, a contradiction. Hence, $x_2x_3 \notin E(G)$. By symmetry, $x_1x_3 \notin E(G)$.

If $x_2x_5$ is an edge, then the graph $F$ shown in Figure~\ref{fig:claim4sub}(c) is a spanning subgraph of $G$. In this case, the set $S = \{x_1,x_5,y_3,y_6\}$ is a dominating set of $F$ (see Figure~\ref{fig:claim4sub}(c)), and so $\gamma \le \gamma(F) = 4$, a contradiction. Hence, $x_2x_5 \notin E(G)$. By symmetry, $x_1x_5 \notin E(G)$.

Renaming $x_1$ and $x_2$ if necessary, we may assume that $x_1x_6$ and $x_2x_4$ are edges. The remaining edge in $G[X]$ is therefore the edge $x_3x_5$. Thus, the graph $G$ is determined.  Renaming vertices if necessary, we may assume that $G = G_{14.2}$, where $G_{14.2}$  is the graph shown in Figure~\ref{fig:claim4sub}(d). We note that $\gamma = 5$. In this case, the set $S = \{y_1,y_3,y_5,v_2\}$ satisfies $\dom_G(S) = 13$ (the vertex $y_4$ represented by the square in Figure~\ref{fig:claim4sub}(d) is the only vertex not dominated by $S$) and $|S| = 4$, implying that $\pd_{\alpha}(G) \le |S| = 4$. This completes the proof of Claim~\ref{c:claim4}.~\smallqed

\medskip
We now return to the proof of Theorem~\ref{t:cubicn14} one final time. As observed earlier, there are four possibilities for the graph $H$, namely $H = 2C_6$ or $H = 3C_4$ or $H = C_4 \cup C_8$ or $H = C_{12}$. By Claim~\ref{c:claim1}, $H \ne 2C_6$. By Claim~\ref{c:claim2}, $H \ne C_{12}$. By Claim~\ref{c:claim3},  if $H = 3C_{4}$, then $\pd_{\alpha}(G) \le 4$. By Claim~\ref{c:claim4}, if $H =  C_4 \cup C_8$, then $\pd_{\alpha}(G) \le 4$. This completes the proof of Theorem~\ref{t:cubicn14}.~\QED

\section{Partial domination in supercubic graphs}
\label{sec:delta-ge-3}

We start this section by proving a useful lemma. We first present a key lemma, which allows us to grow a given set of vertices to a larger set that dominates more vertices. Recall that we refer to graphs $G$ with $\delta(G) \ge 3$ as supercubic graphs.

\begin{lemma} \label{lem:U_S}
Let $k$ be a positive integer and $G$ a supercubic graph of order $n$. If $S \subseteq V(G)$, $U_S= V(G)\setminus N_G[S]$, and
	$$4|U_S| > k(n-|S|)\,,$$
	then there exists a vertex in $\partial(S) \cup U_S$ that dominates at least $k+1$ vertices from $U_S$.
\end{lemma}
\proof Consider the `useful' vertex pairs $(x,y)$ such that $y \in U_S$ and $x$ dominates $y$ (allowing $x = y$). Denote by $p$ the number of useful pairs. As all vertices from $U_S$ can be dominated by itself or one of its at least three neighbors, $p \geq 4|U_S|$. Since $y\in U_S = V(G)\setminus N_G[S]$, we have $N_G[y] \cap S = \emptyset$. 
It follows that $x \in \partial(S) \cup U_S$.

To prove the statement, we suppose that there is no vertex in $G$ which dominates more then $k$ vertices from $U_S$. Equivalently, every vertex $x \in \partial(S) \cup U_S$ belongs to at most $k$ different useful pairs $(x,y)$ (s.t.\ $x$ is the first entry). We then conclude
$$ k(|\partial(S)| + |U_S|) = k(n-|S|) \geq p \geq  4 |U_S|$$
that contradicts the given condition and therefore proves the statement.~\QED

\medskip
We are now in a position to prove that the $\frac{7}{8}$-partial domination number of a cubic graph $G$ is at most one-third of the order of $G$. In fact, our result states that the same is true for every supercubic graph.

\begin{theorem}
\label{t:onethird}
If $G$ is a supercubic graph of order~$n$, then
$$\pd_{\frac{7}{8}}(G) \le \frac{1}{3}n\,.$$
\end{theorem}
\proof First suppose that $G$ is the disjoint union of the components $G_1, \dots , G_k$. It was already observed in~\cite{Das-2019} that $\pd_\alpha (G) \leq \sum_{i=1}^{k} \pd_{\alpha} (G_i)$ holds for each $\alpha$. Therefore, it suffices to prove the statement for connected graphs.

Let $G$ be a connected graph of order~$n$ and of minimum degree $\delta(G) \geq 3$. Let $\alpha = \frac{7}{8}$ and  $\gamma = \gamma(G)$.
We proceed further with two claims.

\begin{claim}
	\label{c:claim12}
	If  $n \le 14$, then $\pd_{\alpha}(G) \le \frac{1}{3}n$.
\end{claim}
\proof  By Theorem~\ref{thm:Reed},
$\gamma \leq \left\lfloor \frac{3}{8}n \right\rfloor$
holds, so $\left\lfloor \frac{3}{8}n \right\rfloor$ vertices are enough to dominate the entire vertex set. Since $\left\lfloor \frac{3}{8}n \right\rfloor = \left\lfloor \frac{1}{3}n \right\rfloor$ holds whenever $n \leq 14 $ and $n \notin \{8,11,14\}$, it suffices  to consider graphs of order $8$, $11$ and $14$.

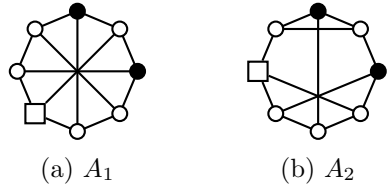
\begin{figure}[ht!]
\begin{center}
\begin{tikzpicture}[scale=.8,style=thick,x=0.8cm,y=0.8cm]
\def\vr{2.5pt} %
\def\vr{3.5pt}
\path (0.00,1.25) coordinate (u1);
\node[rectangle,draw] (u2) at (0.37,0.37) {};
\path (0.37,2.13) coordinate (u8);
\path (1.25,0.00) coordinate (u3);
\path (1.25,2.50) coordinate (u7);
\path (2.13,0.37) coordinate (u4);
\path (2.13,2.13) coordinate (u6);
\path (2.50,1.25) coordinate (u5);
%
\draw (u1)--(u2)--(u3)--(u4)--(u5)--(u6)--(u7)--(u8)--(u1);
\draw (u1)--(u5);
\draw (u2)--(u6);
\draw (u3)--(u7);
\draw (u4)--(u8);
\draw (u1) [fill=white] circle (\vr);
\draw (u3) [fill=white] circle (\vr);
\draw (u4) [fill=white] circle (\vr);
\draw (u5) [fill=black] circle (\vr);
\draw (u6) [fill=white] circle (\vr);
\draw (u7) [fill=black] circle (\vr);
\draw (u8) [fill=white] circle (\vr);
\draw (1.25,-0.8) node {{\small (a) $A_1$}};
\node[rectangle,draw] (v1) at (5.00,1.25) {};
\path (5.37,0.37) coordinate (v2);
\path (5.37,2.13) coordinate (v8);
\path (6.25,0.00) coordinate (v3);
\path (6.25,2.50) coordinate (v7);
\path (7.13,0.37) coordinate (v4);
\path (7.13,2.13) coordinate (v6);
\path (7.50,1.25) coordinate (v5);
%
\draw (v1)--(v2)--(v3)--(v4)--(v5)--(v6)--(v7)--(v8)--(v1);
\draw (v6)--(v8);
\draw (v3)--(v7);
\draw (v1)--(v4);
\draw (v2)--(v5);
%
\draw (v2) [fill=white] circle (\vr);
\draw (v3) [fill=white] circle (\vr);
\draw (v4) [fill=white] circle (\vr);
\draw (v5) [fill=black] circle (\vr);
\draw (v6) [fill=white] circle (\vr);
\draw (v7) [fill=black] circle (\vr);
\draw (v8) [fill=white] circle (\vr);
\draw (6.25,-0.8) node {{\small (b) $A_2$}};
\end{tikzpicture}
\caption{$\frac{7}{8}$-partial dominating sets in $A_1$ and $A_2$. In each, the vertex represented by the square is the only vertex not dominated by the two shaded vertices. }
\label{f:cubic8b}
\end{center}
\end{figure}

Suppose first that $G$ is cubic. Then only $n\in \{8,14\}$ must be considered. If $n = 14$, then by Theorem~\ref{t:cubicn14} we have $\pd_{\alpha}(G) \le \frac{1}{3}n$. Hence we may assume that $n = 8$. If $G$ is isomorphic to $A_1$ or $A_2$, then as illustrated in Figure~\ref{f:cubic8b}, there exists a set $S$ of two (shaded) vertices in $G$ that dominates seven vertices. For any other cubic graph $G$ of order $8$ we have $\gamma(G) \le 2$ by Theorem~\ref{thm:KoSt}. Hence $\pd_{\alpha}(G) \le 2 = \frac{1}{4}n < \frac{1}{3}n$ for each cubic graph $G$ of order $8$.

Assume in the rest that $G$ is supercubic but not cubic. Hence there exists a vertex $u$ of degree at least $4$.

If $n=8$, a vertex $u$ of maximum degree dominates $|N_G[u]|=\dom_G(\{u\}) \ge 5$ vertices. If $\dom_G(\{u\}) \ge 6$,  then any undominated vertex $u'\notin N_G[u]$ can be added to the set and we have $\dom_G(\{u,u' \})\ge 7$. If $\dom_G(\{u\}) =5$, we apply Lemma~\ref{lem:U_S} with $k=1$ and $S=\{u\}$. Since $4|U_S| =4 \times 3 > 8-1 $, there exists a vertex $u'$ such that $\dom_G(\{u,u'\}) \ge 7$. In both cases, $\dom_G(\{u,u'\}) \ge 7$ implies $\pd_\alpha (G) \leq 2 = \left\lfloor \frac{1}{3}n \right\rfloor$.

If $n=11$, we want to prove that there are three vertices $v_1$, $v_2$, $v_3$ that dominate at least $10$ vertices in $G$. Then, $\pd_\alpha (G) \leq 3 = \left\lfloor \frac{1}{3}n \right\rfloor$ will follow. Let $v_1$ be a vertex of maximum degree. We have $\dom_G(\{v_1 \}) \ge 5$. If $\dom_G(\{v_1 \})=5$ then, for $S=\{ v_1\}$, the inequality $4 |U_S|= 24 > 2(n-|S|)=20$ holds and Lemma~\ref{lem:U_S} implies the existence of a vertex $v_2$ that dominates at least three vertices from $U_S$. It follows that $\dom_G(\{v_1, v_2\}) \ge 8$.
If $\dom_G(\{v_1 \})=6$ then, by setting $S=\{ v_1\}$, we get $4 |U_S|= 20 > n-|S|=10$ that shows, by  Lemma~\ref{lem:U_S}, the existence of a vertex $v_2$ which dominates at least two new vertices. Again, we have that $\dom_G(\{v_1, v_2\}) \ge 8$.
If $\dom_G(\{v_1 \})\ge 7$, then $v_2$ can be chosen as an arbitrary undominated vertex and  $\dom_G(\{v_1, v_2\}) \ge 8$ is achieved. For the choice of the last vertex, we consider two cases. If  $\dom_G(\{v_1, v_2\}) = 8$, the set $S=\{v_1, v_2\}$ satisfies the condition in Lemma~\ref{lem:U_S} with $k=1$ and the existence of a vertex $v_3$ which dominates at least two vertices from $U_S$ follows. It means  $\dom_G(\{v_1, v_2, v_3\}) \ge 10$ as required. If $\dom_G(\{v_1, v_2\}) \ge 9$, then any undominated vertex can be chosen as $v_3$ and we have $\dom_G(\{v_1, v_2, v_3\}) \ge 10$ again.

If $n=14$, we want to prove that there exist four vertices $v_1$, $v_2$, $v_3$, $v_4$ which together dominate at least $13$ vertices. Then, $\pd_\alpha (G) \leq 4 = \left\lfloor \frac{1}{3}n \right\rfloor$ will follow. Let $v_1$ be a vertex of maximum degree. If $\dom_G(\{v_1 \})=5$ and $N(v_1)$ is a dominating set in $G$, then  $\dom_G(N(v_1)) =14$ and we are ready. In the other case,  $\dom_G(\{v_1 \})=5$ and $N(v_1)$ is not a dominating set in $G$. Then there exists a vertex $v_2$ such that $\{v_1,v_2\}$ is a packing in $G$. If $v_2$ is a vertex of degree $3$, then $\dom_G(\{v_1,v_2 \} ) = 9$ and for $S=\{v_1, v_2 \}$ and $k=1$, the condition $4 \times 5 > 14- 2$ holds and  Lemma~\ref{lem:U_S} implies the existence of a vertex $v_3$ with  $\dom_G(\{v_1,v_2, v_3 \} ) \ge 11$. If $v_2$ is a vertex of degree at least $4$, then $\dom_G(\{v_1,v_2 \} ) \ge 10$ and $\dom_G(\{v_1,v_2, v_3 \} ) \ge 11$ can be easily achieved. For the choice of the last vertex, we consider two further subcases. If $\dom_G(\{v_1,v_2, v_3 \} ) = 11$, we have $4 \times 3 > 14-3$ and Lemma~\ref{lem:U_S} implies the existence of a vertex $v_4$ with $\dom_G(\{v_1,v_2, v_3, v_4 \} ) \ge 13$. If $\dom_G(\{v_1,v_2, v_3 \} ) \ge 12$ and there are undominated vertices, then we may choose such a vertex $v_4$ and get $\dom_G(\{v_1,v_2, v_3, v_4 \} ) \ge 13$.~\smallqed

\bigskip
By Claim~\ref{c:claim12}, we may assume that $n \ge 15$, for otherwise the desired result follows. Let $D=\{v_1,v_2,\ldots,v_{\gamma}\}$ be a $\gamma$-set of $G$ satisfying the Bollob\'{a}s-Cockayne Lemma~\ref{lem:BoCo}, and so $\epn[v,D] \ne \emptyset$ for every vertex $v \in D$. By Theorem~\ref{thm:Reed}, we have $\gamma \le \lfloor \frac{3}{8}n \rfloor$. If $\gamma \le \frac{1}{3}n$, then the set $D$ is certainly an $\alpha$-partial dominating set of $G$ of cardinality at most $\frac{1}{3}n$. Thus in this case, $\pd_{\alpha}(G) \le |D| \le \frac{1}{3}n$. Hence we may assume that $\gamma > \frac{1}{3}n$, for otherwise the desired result is immediate.

Using the vertices $v_1, \dots , v_\gamma$ from $D$, let $(V_1,V_2,\ldots,V_{\gamma})$ be a partition of the vertex set $V(G)$ such that for all $i \in [\gamma]$, the following properties hold: (i) $v_i \in V_i$, (ii) $\epn[v_i,D] \subset V_i$, and (iii) $V_i \subseteq N_G[v_i]$. As observed earlier,  $|\epn[v_i,D]| \ge 1$, and so $|V_i| \ge |\{v_i\}| + |\epn[v_i,D]| \ge 2$ for all $i \in [\gamma]$. Renaming the vertices $v_1, v_2, \ldots, v_{\gamma}$ if necessary, we may assume that $|V_i| \ge |V_{i+1}|$ for all $i \in [\gamma - 1]$, that is,
\begin{equation}
\label{Eq2}
|V_1| \ge |V_2| \ge \cdots \ge |V_{\gamma}| \ge 2.
\end{equation}

Let $k_1 = \lfloor \frac{1}{3}n \rfloor$ and let $k_2 = \gamma - k_1$. By assumption, $\frac{1}{3}n < \gamma$. By Theorem~\ref{thm:Reed}, $\gamma \le \lfloor \frac{3}{8}n \rfloor$. Hence, $\frac{1}{3}n < \gamma \le  \lfloor \frac{3}{8}n \rfloor$. By definition of $k_1$ and $k_2$ and by our earlier observations and assumptions,
\begin{equation}
\label{Eq3}
1 \le k_2 = \gamma - k_1 \le  \left\lfloor \frac{3}{8}n \right\rfloor - \left\lfloor \frac{1}{3}n \right\rfloor.
\end{equation}

Let $S = \{v_1,v_2,\ldots,v_{k_1}\}$, and so $|S| = k_1 = \lfloor \frac{1}{3}n \rfloor$. Since $(V_1,V_2,\ldots,V_{\gamma})$ is a partition of the vertex set $V(G)$, we note that the number of vertices dominated by the set $S$ is at least the number of vertices in the sets $V_1 \cup \cdots \cup V_{k_1}$, that is,
\begin{equation}
\label{Eq4}
\dom_G(S) \ge \sum_{i=1}^{k_1} |V_i|.
\end{equation}

We proceed further with the following claim.

\begin{claim}
\label{c:claim13}
If $|V_{k_1}| \ge 3$, then $\pd_{\alpha}(G) \le \frac{1}{3}n$.
\end{claim}
\proof Suppose that $|V_{k_1}| \ge 3$. In this case, by Inequalities~(\ref{Eq2}) and~(\ref{Eq4}), and by our assumption that $n \ge 15$, we infer that
\begin{equation}
\label{Eq5}
\dom_G(S) \ge 3k_1 = 3\left\lfloor \frac{1}{3}n \right\rfloor \ge \left\lceil \frac{7}{8}n \right\rceil. \2
\end{equation}

By Inequality~(\ref{Eq5}), we have $\dom_G(S) \ge \frac{7}{8}n$, implying that the set $S$ is an $\alpha$-partial dominating set of $G$, and so $\pd_{\alpha}(G) \le |S| \le \frac{1}{3}n$, yielding the desired result.~\smallqed

\medskip
By Claim~\ref{c:claim13}, we may assume that $|V_{k_1}| = 2$, for otherwise the desired result follows. With this assumption and by inequality~(\ref{Eq2}), we note that $|V_i| = 2$ for all $i \ge k_1$. Hence by Inequality~(\ref{Eq3}), we have
\begin{equation}
\label{Eq6}
\sum_{i=k_1 + 1}^{k_2} |V_i| = 2k_2 \le 2\left(\left\lfloor \frac{3}{8}n \right\rfloor - \left\lfloor \frac{1}{3}n \right\rfloor \right).
\end{equation}

Thus, by inequalities~(\ref{Eq4}) and~(\ref{Eq6}), and by our assumption that $n \ge 15$, we infer that
\begin{equation}
\label{Eq7}
\dom_G(S) \ge \sum_{i=1}^{k_1} |V_i| = n - \sum_{i=k_1 + 1}^{k_2} |V_i| \ge n - 2\left(\left\lfloor \frac{3}{8}n \right\rfloor - \left\lfloor \frac{1}{3}n \right\rfloor \right) \ge \left\lceil \frac{7}{8}n \right\rceil.
\end{equation}

By Inequality~(\ref{Eq7}), we have $\dom_G(S) \ge \frac{7}{8}n$, implying that the set $S$ is an $\alpha$-partial dominating set of $G$, and so $\pd_{\alpha}(G) \le |S| \le \frac{1}{3}n$, yielding the desired result.~\QED

\medskip
The bound in Theorem~\ref{t:onethird} is best possible in the sense that if $\alpha > \frac{7}{8}$ and $G$ is $A_1$ or $A_2$ (see Figure~\ref{f:cubic8}), then in this case $\lceil \alpha \times n \rceil = 8 = n$, and at least three vertices are needed to dominate all vertices of $G$. Thus in this example, $\pd_{\alpha}(G) = 3 = \frac{3}{8}n > \frac{1}{3}n$. The same is true if every component of $G$ is isomorphic to $A_1$ or $A_2$. Hence the value for $\alpha$ in the statement of Theorem~\ref{t:onethird} cannot be increased in general in order to guarantee that the $\alpha$-partial domination number of a connected cubic graph is at most one-third its order.

However if the connected cubic graph $G$ has sufficiently large order~$n$, then we can improve the value $\alpha = \frac{7}{8}$ given in Theorem~\ref{t:onethird} to a larger value of $\alpha$. For example, if $n \ge 28$, then $\alpha = \frac{13}{14}$ suffices, as the following result shows.

\begin{theorem}
\label{t:onethird-improved}
If $G$ is a connected cubic graph of order $n \ge 28$, then
$$\pd_{\frac{13}{14}}(G) \le \frac{1}{3}n\,.$$
\end{theorem}
\proof Let $G$ be a connected cubic graph of order~$n \ge 28$ and let $\alpha = \frac{13}{14}$. We adopt exactly our notation from the proof of Theorem~\ref{t:onethird}. In particular, $D$ is a $\gamma$-set of $G$ satisfying Lemma~\ref{lem:BoCo}. As before, by Theorem~\ref{thm:KoSt} we have $\gamma \le \lfloor \frac{5}{14}n \rfloor$. If $\gamma \le \frac{1}{3}n$, then $\dom_G(D) = n$. Hence we may assume that $\gamma > \frac{1}{3}n$, for otherwise the desired result is immediate. Let $k_1$ and $k_2$ be defined exactly as in the proof Theorem~\ref{t:onethird}. If $|V_{k_1}| \ge 3$, then
\begin{equation}
\label{eq:first}
\dom_G(S) \ge 3k_1 = 3\left\lfloor \frac{1}{3}n \right\rfloor \ge \left\lceil \frac{13}{14}n \right\rceil, \2
\end{equation}
implying that the set $S$ is an $\alpha$-partial dominating set of $G$. Thus, $\pd_{\alpha}(G) \le |S| \le \frac{1}{3}n$, yielding the desired result. Hence we may assume that $|V_{k_1}| = 2$. With this assumption, we note that $|V_i| = 2$ for all $i \ge k_1$. Thus proceeding exactly as before, we yield the inequality chain where recall that by supposition we have $n \ge 28$.
\begin{equation}
\label{eq:second}
\dom_G(S) \ge \sum_{i=1}^{k_1} |V_i| = n - \sum_{i=k_1 + 1}^{k_2} |V_i| \ge n - 2\left(\left\lfloor \frac{5}{14}n \right\rfloor - \left\lfloor \frac{1}{3}n \right\rfloor \right) \ge \left\lceil \frac{13}{14}n \right\rceil.
\end{equation}
Once again, $\dom_G(S) \ge \frac{13}{14}n$, implying that the set $S$ is an $\alpha$-partial dominating set of $G$. Thus, $\pd_{\alpha}(G) \le |S| \le \frac{1}{3}n$.~\QED

Note that in the proof of Theorem~\ref{t:onethird-improved} we used the inequality $\gamma \le \lfloor \frac{5}{14}n \rfloor$ which holds for every connected cubic graph of order at least~$10$. Hence we cannot avoid the assumption that $G$ is connected. On the other hand, $\gamma \le \lfloor \frac{3}{8}n \rfloor$ holds for every supercubic graph, and we have the following result.

\begin{theorem}
\label{t:onethird-improved-delta-at-least-3}
If $G$ is a supercubic graph of order~$n \ge 60$, then
$$\pd_{\frac{9}{10}}(G) \le \frac{1}{3}n\,.$$
\end{theorem}
\proof
We can proceed along the same lines as in the proof of Theorem~\ref{t:onethird-improved}. The only difference is that now we cannot apply Theorem~\ref{thm:KoSt}, instead we apply Theorem~\ref{thm:Reed}. Then~\eqref{eq:first} rewrites as
\begin{equation}
\label{eq:first'}
\dom_G(S) \ge 3k_1 = 3\left\lfloor \frac{1}{3}n \right\rfloor \ge \left\lceil \frac{9}{10}n \right\rceil, \2
\end{equation}
which holds for $n\ge 18$, while~\eqref{eq:second} rewrites as
\begin{equation}
\label{eq:second'}
\dom_G(S) \ge \sum_{i=1}^{k_1} |V_i| = n - \sum_{i=k_1 + 1}^{k_2} |V_i| \ge n - 2\left(\left\lfloor \frac{3}{8}n \right\rfloor - \left\lfloor \frac{1}{3}n \right\rfloor \right) \ge \left\lceil \frac{9}{10}n \right\rceil,
\end{equation}
which holds for $n\ge 60$. The conclusion follows.~\QED

\section{Closing remarks}
\label{S:closing}

As a consequence of Theorems~\ref{thm:Reed} and~\ref{thm:KoSt}, we have the following result which characterizes the connected cubic graphs $G$ of order~$n$ satisfying $\gamma(G) = \frac{3}{8}n$.

\begin{corollary}{\rm (\cite{KoStocker-09,Re-96})}
If $G$ is a connected cubic graph of order $n$, then $\gamma(G) \le \frac{3}{8}n$, with equality if and only if $G$ is one of the two graphs $A_1$ and $A_2$ shown in Figure~\ref{f:cubic8}.
\end{corollary}

A natural problem is to characterize the graphs that achieve equality in the Kostochka-Stocker Theorem~\ref{thm:KoSt}, that is, to characterize the connected cubic graphs $G$ of order~$n$ satisfying $\gamma(G) = \frac{5}{14}n$. Necessarily, for such graphs we have $n = 14k$ for some $k \ge 1$.

\begin{figure}[ht!]
\begin{center}
\begin{tikzpicture}[scale=0.65,style=thick]
\tikzstyle{every node}=[draw=none,fill=none]
\def\vr{3pt} 

\begin{scope}[yshift = 0cm, xshift = 0cm]
\path (0,0) coordinate (x1);
\path (1.5,0) coordinate (x2);
\path (3,0) coordinate (x3);
\path (4.5,0) coordinate (x4);
\path (6,0) coordinate (x5);
\path (7.5,0) coordinate (x6);
\path (9,0) coordinate (x7);
\path (0,2) coordinate (y1);
\path (1.5,2) coordinate (y2);
\path (3,2) coordinate (y3);
\path (4.5,2) coordinate (y4);
\path (6,2) coordinate (y5);
\path (7.5,2) coordinate (y6);
\path (9,2) coordinate (y7);
\draw (x1) --(x2);
\draw (x6)--(x7);
\draw (y1) --(y2)--(y3)--(y4)--(y5)--(y6)--(y7);
\draw (y1) .. controls (1,3) and (8,3) .. (y7);
\draw (x1) .. controls (0.5,-1) and (2.5,-1) .. (x3);
\draw (x2) .. controls (2.0,-1) and (4.0,-1) .. (x4);
\draw (x3) .. controls (3.5,-1) and (5.5,-1) .. (x5);
\draw (x4) .. controls (5.0,-1) and (7.0,-1) .. (x6);
\draw (x5) .. controls (6.5,-1) and (8.5,-1) .. (x7);

\foreach \i in {1,...,7}
\draw (x\i)--(y\i);
\foreach \i in {1,...,7}
{
\draw (x\i)  [fill=black] circle (\vr);
\draw (y\i)  [fill=black] circle (\vr);
}
\draw (4.5,-1.5) node {{\small (a)} $G_{14.3}$};
\end{scope}

\begin{scope}[yshift = 0.00cm, xshift = 12cm]
\path (0,0) coordinate (x1);
\path (1.5,0) coordinate (x2);
\path (3,0) coordinate (x3);
\path (4.5,0) coordinate (x4);
\path (6,0) coordinate (x5);
\path (7.5,0) coordinate (x6);
\path (9,0) coordinate (x7);
\path (0,2) coordinate (y1);
\path (1.5,2) coordinate (y2);
\path (3,2) coordinate (y3);
\path (4.5,2) coordinate (y4);
\path (6,2) coordinate (y5);
\path (7.5,2) coordinate (y6);
\path (9,2) coordinate (y7);
\draw (x1) --(x2)--(x3)--(x4)--(x5)--(x6)--(x7);
\draw (y1) --(y2)--(y3)--(y4)--(y5)--(y6)--(y7);
\draw (x1) --(y1);
\draw (x2) --(y2);
\draw (x3) --(y4);
\draw (x4) --(y3);
\draw (x5) --(y5);
\draw (x6) --(y7);
\draw (x7) --(y6);
\draw (x1) .. controls (1,-1.0) and (8,-1.0) .. (x7);
\draw (y1) .. controls (1,3.0)and (8,3.0) .. (y7);

\foreach \i in {1,...,7}
{
\draw (x\i)  [fill=black] circle (\vr);
\draw (y\i)  [fill=black] circle (\vr);
}
\draw (4.5,-1.75) node {{\small (b)} $G_{14.2}$};
\end{scope}

\begin{scope}[yshift = -6cm, xshift = -0.5cm]
\path (0,0) coordinate (x1);
\path (1.5,0) coordinate (x2);
\path (3,0) coordinate (x3);
\path (4.5,0) coordinate (x4);
\path (6,0) coordinate (x5);
\path (7.5,0) coordinate (x6);
\path (0,2) coordinate (y1);
\path (1.5,2) coordinate (y2);
\path (3,2) coordinate (y3);
\path (4.5,2) coordinate (y4);
\path (6,2) coordinate (y5);
\path (7.5,2) coordinate (y6);
\path (9,1) coordinate (w);
\path (10.5,1) coordinate (t);
\draw (x1) --(x2)--(x3)--(x4)--(x5)--(x6);
\draw (y1) --(y2)--(y3)--(y4)--(y5)--(y6);
\draw (x1) --(y2);
\draw (x2) --(y1);
\draw (x3) --(y3);
\draw (x4) --(y4);
\draw (x5) --(y6);
\draw (x6) --(y5);
\draw (x6) --(w);
\draw (y6) --(w);
\draw (w) --(t);
\draw (y1) .. controls (1,3.5) and (8,3.5) .. (t);
\draw (x1) .. controls (1,-1.5)and (8,-1.5) .. (t);
\foreach \i in {1,...,6}
{
\draw (x\i)  [fill=black] circle (\vr);
\draw (y\i)  [fill=black] circle (\vr);
}
\draw (w)  [fill=black] circle (\vr);
\draw (t)  [fill=black] circle (\vr);
\draw (4.5,-1.75) node {{\small (c) $G_{14.1}$}};
\draw (16.5,-1.75) node {{\small (d) $P(7,2)$}};

\end{scope}

\begin{scope}[yshift = -5cm, xshift = 16cm]
\tikzstyle{every node}=[circle, draw, fill=black!100, inner sep=0pt,minimum width=.135cm]
  \draw \foreach \x in {38.57,90,...,347.15}
  {
    (\x:2) node{} -- (\x+51.43:2)
    (\x:2) -- (\x:1) node{}
    (\x:1) -- (\x+102.86:1)
  };
\end{scope}

\end{tikzpicture}
\end{center}
\vskip -0.5 cm
\caption{The four connected cubic graphs $G$ of order~$n = 14$ satisfying $\gamma(G) = 5$}
\label{fig:X}
\end{figure}
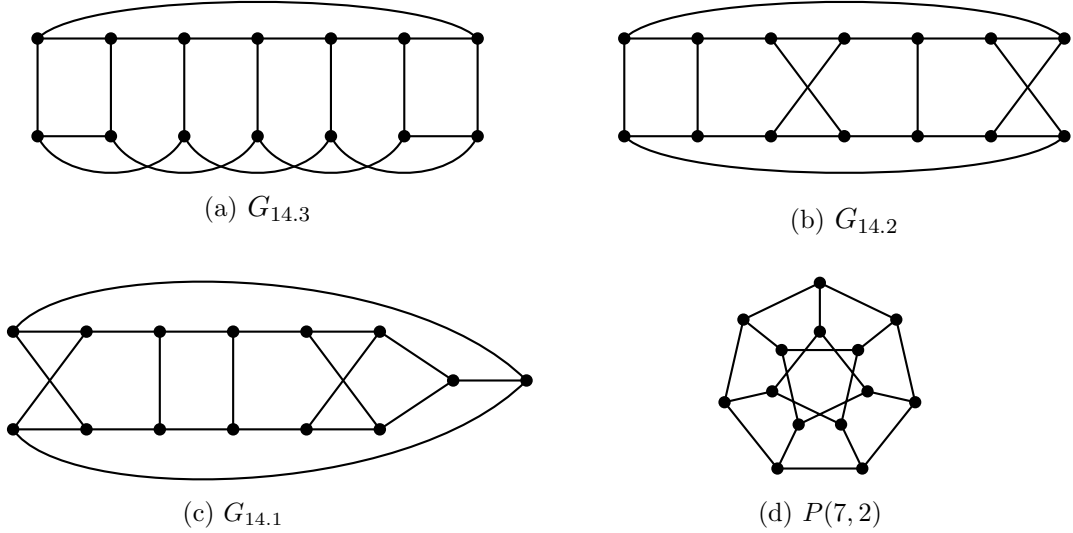

We show next that there are exactly four such graphs of order~$n=14$. We remark that the proof of Theorem~\ref{t:cubicn14} gave rise to two connected cubic graphs of order~$n$ satisfying $\gamma(G) = 5 = \frac{5}{14}n$, namely the graphs $G_{14.1}$ and $G_{14.2}$ shown in Figures~\ref{fig:claim3}(d) and~\ref{fig:claim4sub}(d), respectively. (These two graphs are redrawn in Figure~\ref{fig:X}(c) and~\ref{fig:X}(b), respectively.) With a bit more work, one can readily establish two additional such graphs.

In the second paragraph of the proof of Theorem~\ref{t:cubicn14}, we consider the case when $\rho(G) = 3$. In this case, adding a vertex at distance~$3$ to a maximum packing immediately yielded an $\frac{7}{8}$-partial dominating set of $G$ of cardinality~$4$, and therefore we assumed that $\rho(G) = 2$. However a more detailed analysis of the case when $\rho(G) = 3$ yields the generalized Petersen graph $P(7,2)$ shown in Figure~\ref{fig:X}(d).

In the fourth paragraph of the proof of Theorem~\ref{t:cubicn14}, we considered the case when the set $X = V(G) \setminus N_G[P]$ contains a vertex adjacent to two other vertices in $X$. Since this case immediately yielded an $\frac{7}{8}$-partial dominating set of $G$ of cardinality~$4$, we therefore assumed that this case does not occur. However a more detailed analysis of the case when a vertex in $X$ has two neighbors in $X$ yields the graph $G_{14.3}$ shown in Figure~\ref{fig:X}(a). The proof details giving rise to these two additional graphs, namely $P(7,2)$ and $G_{14.3}$, are similar to our proof of Theorem~\ref{t:cubicn14}, and are not given here. Moreover, the result was also verified by computer.

\begin{theorem}
\label{t:largedom}
If $G$ is a connected cubic graph of order~$n = 14$ satisfying $\gamma(G) = 5 = \frac{5}{14}n$, then $G \in \{G_{14.1},G_{14.2},G_{14.3},P(7,2)\}$.
\end{theorem}

It is not known if the $\frac{5}{14}$-upper bound on the domination number of a connected cubic graph of order~$n$ given by Kostochka and Stocker~\cite{KoStocker-09} is achievable when $n$ is large. We pose the following conjecture.

\begin{conj}
\label{c:largedom}
If $G$ is a connected cubic graph of order~$n$ satisfying $\gamma(G) = \frac{5}{14}n$, then $G \in \{G_{14.1},G_{14.2},G_{14.3},P(7,2)\}$.
\end{conj}

The authors in~\cite{KoStocker-09} remark that the bound $\gamma(G) \le \lfloor \frac{5}{14}n \rfloor$ for a connected cubic graph of order~$n \ge 14$ is achievable for $n \in \{14,16,18\}$. It would be interesting to find graphs of orders $n \ge 20$ that achieve equality in this bound. Natural candidates are the generalized Petersen graphs $P(p,2)$ of order~$n = 2p$ whose domination numbers are known  (see,~\cite{FuYuJi-09}).

\begin{theorem}{\rm (\cite{FuYuJi-09})}
\label{thm:GPeteresen}
$\gamma(P(p,2)) = p - \left\lfloor \frac{p}{5} \right\rfloor - \left\lfloor \frac{p+2}{5} \right\rfloor$ for all $p \ge 3$.
\end{theorem}

For $p \in \{3,5,6,7,8,9,11,12\}$, we have $p - \lfloor \frac{p}{5} \rfloor - \lfloor \frac{p+2}{5} \rfloor = \lfloor \frac{5}{7}p \rfloor$. Hence as a consequence of Theorem~\ref{thm:GPeteresen}, we have the following result.

\begin{corollary}
\label{cor:GPeteresen}
For $n \in \{6,10,12,14,16,18,22,24\}$, there exist connected cubic graphs $G$ of order~$n$ satisfying $\gamma(G) = \lfloor \frac{5}{14}n \rfloor$.
\end{corollary}

As far as we are aware, $P(12,2)$ is the largest currently known connected cubic graph of order~$n$ satisfying $\gamma(G) = \lfloor \frac{5}{14}n \rfloor$. In addition, $\gamma(P(12,2)) = 8 = \frac{1}{3}n$. We close with the following question, for which we suspect the answer is no. 

\begin{quest}
\label{q:largen}
Are there infinitely many connected cubic graphs $G$ of order~$n$ satisfying $\gamma(G) = \lfloor \frac{5}{14}n \rfloor$?
\end{quest}

\section*{Acknowledgements}

We are grateful to Gregor Rus for computer verification of Theorem~\ref{t:cubicn14}. This work was supported by the Slovenian Research Agency (ARRS) under the grants P1-0297, J1-2452, and N1-0285. Research of the second author was supported in part by the University of Johannesburg and the South African National Research Foundation.


\begin{thebibliography}{99}

\bibitem{BoCo-79}
B.~Bollob{\'a}s, E.~J. Cockayne,
Graph-theoretic parameters concerning domination, independence, and irredundance,
J.\ Graph Theory 3 (1979) 241--249.

\bibitem{Borg-2020}
P.~Borg, P.~Kaemawichanurat,
Partial domination of maximal outerplanar graphs,
Discrete Appl.\ Math.\ 283 (2020) 306--314.


\bibitem{CaIs-21}
J.~M.~C. Cabulao, R.~T.~Isla,
On connected partial domination in graphs,
Eur. J. Pure Appl. Math. 14 (2021) 1490--1506.



\bibitem{Case-2017}
B.~M.~Case, S.~T.~Hedetniemi, R.~C.~Laskar, D.~J.~Lipman,
Partial domination in graphs,
Congr.\ Numer.\ 228 (2017) 85--95.

\bibitem{Das-2019}
A.~Das,
Partial domination in graphs,
Iran.\ J.\ Sci.\ Technol.\ Trans.\ A Sci.\ 43 (2019) 1713--1718.

\bibitem{DaDe-2018}
A.~Das, W.~J.~Desormeaux,
Domination defect in graphs: guarding with fewer guards,
Indian J. Pure Appl. Math. 49 (2018) 349--364.

\bibitem{DaLa-2018}
A.~Das, R.~C.~Laskar, N. Jafari Rad,
On $\alpha$-domination in graphs,
Graphs Combin. 34 (2018) 193--205.



\bibitem{Dunbar-2000}
J.~E.~Dunbar, D.~G.~Hoffman, R.~C.~Laskar, L.~R.~Markus,
$\alpha$-domination,
Discrete Math. 211 (2000) 11--26.

\bibitem{Fa-96}
O.~Favaron,
Signed domination in regular graphs,
Discrete Math.\ 158 (1996) 287--293.

\bibitem{FuYuJi-09}
X.~Fu, Y.~Yang, B.~Jiang,
On the domination number of generalized Petersen graphs $P(n,2)$,
Discrete Math. 309 (2009) 2445--2451.

\bibitem{HaHeHe-20}
T.~W.~Haynes, S.~T.~Hedetniemi, M.~A.~Henning (eds.),
Topics in Domination in Graphs.
Developments in Mathematics 64, Springer, Cham, 2020.

\bibitem{HaHeHe-21}
T.~W.~Haynes, S.~T.~Hedetniemi, M.~A.~Henning (eds.),
Structures of Domination in Graphs.
Developments in Mathematics 66, Springer, Cham, 2021.

\bibitem{HaHeHe-23}
T.~W.~Haynes, S.~T.~Hedetniemi, M.~A.~Henning,
Domination in Graphs: Core Concepts.
Springer Monographs in Mathematics, Springer, Cham, 2023.

\bibitem{Henning-2022}
M.~A.~Henning,
Bounds on domination parameters in graphs: a brief survey,
Discuss.\ Math.\ Graph Theory 42 (2022) 665--708.

\bibitem{RadVolk-2016}
N.~Jafari Rad, L.~Volkmann, Lutz
Edge-removal and edge-addition in $\alpha$-domination,
Graphs Combin. 32 (2016) 1155--1166.

\bibitem{Ke-06}
A.~Kelmans,
Counterexamples to the cubic graph domination conjecture,
http://arxiv.org/pdf/math/0607512.pdf (20 July 2006).

\bibitem{KoSt-05}
A.~V.~Kostochka, B.~Y.~Stodolsky,
On domination in connected cubic graphs,
Discrete Math.\ 304 (2005) 45--50.

\bibitem{KoSt-09}
A.~V.~Kostochka, B.~Y.~Stodolsky,
An upper bound on the domination number of $n$-vertex connected cubic graphs,
Discrete Math.\ 309 (2009) 1142--1162.

\bibitem{KoStocker-09}
A.~V.~Kostochka, C.~Stocker,
A new bound on the domination number of connected cubic graphs,
Sib.\ Elektron.\ Mat.\ Izv.\ 6 (2009) 465--504.

\bibitem{MaIs-10}
R.~D.~Macapodi, R.~Y.~Isla,
Total partial domination in graphs under some binary operations,
Eur. J. Pure Appl. Math. 12 (2019) 1643--1655.

\bibitem{Kral-2012}
D.~Kr\'al, P.~\v{S}koda, J.~Volec,
Domination number of cubic graphs with large girth,
J.\ Graph Theory 69 (2012) 131--142.

\bibitem{Lowenstein-2008}
C.~L\"{o}wenstein, D.~Rautenbach,
Domination in graphs of minimum degree at least two and large girth,
Graphs Combin.\ 24 (2008) 37--46.

\bibitem{PhVa-21}
N.~L.~Philo, J. Varghese Kureethara,
Independent partial domination,
Cubo 23 (2021) 411--421.

\bibitem{Rautenbach-2008}
D.~Rautenbach, B.~Reed,
Domination in cubic graphs of large girth,
Lecture Notes in Comput.\ Sci.\ 4535 (Springer, Berlin, 2008)
186--190.

\bibitem{Re-96}
B.~Reed,
Paths, stars, and the number three,
Combin.\ Probab.\ Comput.\ 5 (1996) 277--295.

\end{thebibliography}
\end{document}